\documentclass{amsart}
\usepackage[alphabetic]{amsrefs}
\usepackage{accents}
\usepackage{graphicx}
\usepackage{amsfonts}
\usepackage{amscd}
\usepackage{amssymb}
\usepackage{xypic}
\xyoption{all}
\usepackage{hyperref}
\usepackage{url}
\usepackage{tikz}
\usepackage{calligra}
\usepackage{extpfeil}
\usetikzlibrary{shapes}
\usetikzlibrary{decorations.markings}
\usepackage{tikz-cd}
\usetikzlibrary{arrows}
\usepackage{multirow}
\usepackage[savepos]{zref}
\newtheorem{thm}{Theorem}[section]
\newtheorem{corollary}[thm]{Corollary}
\newtheorem{lemma}[thm]{Lemma}
\newtheorem{proposition}[thm]{Proposition}

\newtheorem{example}[thm]{Example}

\theoremstyle{definition}
\newtheorem{definition}[thm]{Definition}
\theoremstyle{remark}
\newtheorem{remark}[thm]{Remark}
\numberwithin{equation}{section}
\numberwithin{figure}{section}

\usepackage[mathscr]{euscript}
\usepackage[T2A,T1]{fontenc}

  {\begin{description}%
    \setlength{\itemsep}{5pt}%
    \setlength{\parskip}{5pt}}%
  {\end{description}}

\newcommand{\Sha}{\mbox{\usefont{T2A}{\rmdefault}{m}{n}\CYRSHCH}}

\newcommand{\form}{\upsilon}

\DeclareMathOperator{\Fr}{{Fr}}

\DeclareMathOperator{\unit}{u}
\DeclareMathOperator{\fin}{fin}
\DeclareMathOperator{\inc}{inc}

\DeclareMathOperator{\mult}{m}

\DeclareMathOperator{\Pic}{Pic}

\DeclareMathOperator{\Comm}{Comm}

\DeclareMathOperator{\Ind}{Ind}

\DeclareMathOperator{\zind}{Zind}
\DeclareMathOperator{\cind}{Cind}
\DeclareMathOperator{\Hom}{Hom}
\DeclareMathOperator{\Sym}{Sym}

\DeclareMathOperator{\End}{End}

\DeclareMathOperator{\Gal}{Gal}
\DeclareMathOperator{\Span}{Span}
\DeclareMathOperator{\unr}{unr}

\DeclareMathOperator{\Ker}{Ker}

\DeclareMathOperator{\val}{val}

\DeclareMathOperator{\Hilb}{Hilb}
\DeclareMathOperator{\Char}{char}

\DeclareMathOperator{\Spec}{Spec}

\DeclareMathOperator{\Id}{Id}
\DeclareMathOperator{\Res}{Res}

\DeclareMathOperator{\Tr}{Tr}
\DeclareMathOperator{\Kum}{\varkappa}

\newdir{ >}{{}*!/-5pt/@{>}}
\newdir^{ (}{{}*!/-5pt/@^{(}}

\DeclareMathOperator{\et}{\acute{e}t}

\newcommand{\loccit}{loc.~cit.}
\newcommand{\From}{\colon}
\newcommand{\inar}{\ar@{^{(}->}}
\newcommand{\onar}{\ar@{->>}}

\renewcommand{\Im}{\mathrm{Im}}

\newcommand{\defined}[1]{\underline{{#1}}}
\usepackage{accents}
\newlength{\dtildeheight}
\newcommand{\dtilde}[1]{%
    \settoheight{\dtildeheight}{\ensuremath{\tilde{#1}}}%
    \addtolength{\dtildeheight}{-0.15ex}%
    \tilde{\vphantom{\rule{1pt}{\dtildeheight}}%
    \smash{\tilde{#1}}}}

\newcommand{\Irr}{{\boldsymbol{\Pi}}}
\newcommand{\Baer}{\dotplus}

\newcommand{\limdir}{\varinjlim}

\newcommand{\Cat}[1]{ {\mathsf{#1}} }
\newcommand{\Fun}[1]{ {\mathsf{#1}} }

\newcommand{\Lie}[1]{ {\mathfrak{#1}} }

\newcommand{\alg}[1]{\boldsymbol{\mathrm{#1}}}

\newcommand{\sheaf}[1]{{\mathscr{#1}}}
\newcommand{\shom}{\mathscr{H}\mathit{om}}
\newcommand{\ssym}{\mathscr{S}\mathit{ym}}

\newcommand{\amu}{\alg{\mu}}

\newcommand{\Ext}{ {\mbox{Ext}} }

\newcommand{\ZZ}{\mathbb Z}

\newcommand{\RR}{\mathbb R}

\newcommand{\CC}{\mathbb C}

\renewcommand{\AA}{\mathbb A}
\newcommand{\VV}{\mathcal V}
\newcommand{\WW}{\mathcal W}

\newcommand{\UU}{\mathcal U}

\renewcommand{\SS}{\mathcal S}

\newcommand{\hecke}{\mathcal H}

\newcommand{\FF}{\mathbb F}
\newcommand{\OO}{\mathcal{O}}

\newcommand{\Into}{\hookrightarrow}
\newcommand{\Onto}{\twoheadrightarrow}
\newcommand{\To}{\rightarrow}

\newcommand{\isom}{\cong}

\newcommand{\AF}{\mathcal{AF}}

\renewcommand{\th}{\text{th}}
\newcommand{\Pouch}{\mathcal{P}}

\DeclareMathOperator{\loc}{loc}
\DeclareMathOperator{\TD}{Tate}
\newcommand{\inarrow}{\arrow[hook]}
\newcommand{\onarrow}{\arrow[two heads]}
\newcommand{\PD}{\ast}

\makeatletter
\newcommand\@biprod[1]{%
  \vcenter{\hbox{\ooalign{$#1\prod$\cr$#1\coprod$\cr}}}}
\newcommand\biprod{\mathop{\mathpalette\@biprod\relax}\displaylimits}
\makeatother

\newcommand{\defeq}{:=}
\DeclareMathAlphabet{\mathcalligra}{T1}{calligra}{m}{n}
\DeclareMathOperator{\Zar}{Zar}
\DeclareMathOperator{\zar}{zar}

\begin{document}

\title{Covers of tori over local and global fields}%
\author{Martin H. Weissman}%
\date{\today}

\address{Yale-NUS College, 6 College Ave East, \#B1-01, Singapore 138614}
\email{marty.weissman@yale-nus.edu.sg}%


\begin{abstract}
Langlands has described the irreducible admissible representations of $T$, when $T$ is the group of points of an algebraic torus over a local field.  Also,  Langlands described the automorphic representations of $T_\AA$ when $T_\AA$ is the group of adelic points of an algebraic torus over a global field $F$.  

We describe irreducible (in the local setting) and automorphic (in the global setting) $\epsilon$-genuine representations for ``covers'' of tori, also known as ``metaplectic tori,'' which arise from a framework of Brylinski and Deligne.  In particular, our results include a description of spherical Hecke algebras in the local unramified setting, and a global multiplicity estimate for automorphic representations of covers of split tori.  For automorphic representations of covers of split tori, we prove a multiplicity-one theorem.
\end{abstract}

\maketitle

\tableofcontents

\section*{Introduction}

This article studies the representation theory of ``covers'' of tori over local and global fields.  For a torus $\alg{T}$ over a local field $F$, a central extension $\alg{K}_2 \Into \alg{T}' \Onto \alg{T}$, and a suitable positive integer $n$ (the condition is compactly stated as $\# \amu_n(F) = n$), the construction of Brylinski and Deligne \cite[Construction 10.3]{B-D} gives a central extension
$$\mu_n \Into \tilde T \Onto T = \alg{T}(F).$$
For such an extension $\alg{K}_2 \Into \alg{T}' \Onto \alg{T}$ over a global field $F$, their construction \cite[\S 10.4]{B-D} gives a central extension
$$\mu_n \Into \tilde T_\AA \Onto T_\AA,$$
together with a canonical splitting over $T_F$.  An injective character $\epsilon \From \mu_n \Into \CC^\times$ is fixed throughout the paper.  Our previous work \cite{MWTori} discussed the $\epsilon$-genuine irreducible representations of $\tilde T$ in the local case.  Here we go further, classifying irreducible $\epsilon$-genuine admissible representations of $\tilde T$ (in the local case) or $\tilde T_\AA$ (in the global case) by their central character.  

Section \ref{Center} provides a description of $Z(\tilde T)$ (in the local case) and $Z(\tilde T_\AA)$ (in the global case).  This description is precise, in terms of the image of an isogeny (on $F$-points or $\AA$-points), when $\alg{T}$ is split or quasitrivial (i.e., obtained by Weil restriction of scalars from a split torus).  Section \ref{Irreps} describes the irreducible admissible representations of $\tilde T$ (in the local case) in terms of their central character -- this only goes slightly beyond \cite{MWTori}, and is provided mostly to keep the paper self-contained.  We partition the irreducible $\epsilon$-genuine admissible representations of $\tilde T$ into {\em pouches}.  We expect elements of the same pouch to be ``L-indistinguishable,'' i.e., to have the same parameter in a future classification by Weil-parameters.

When $\alg{T}$ is a torus over the ring $\OO$ of integers in a nonarchimedean local field $F$, we have considered extensions $\alg{K}_2 \Into \alg{T}' \Onto \alg{T}$ over $\OO$ in \cite{MWIntegral}.  For such an extension, observations from \cite[\S 10.7]{B-D} demonstrate that the central extension $\mu_n \Into \tilde T \Onto T$ is equipped with a splitting over $T^\circ = \alg{T}(\OO)$.  This enables a definition of {\em unramified} $\epsilon$-genuine representations -- those generated by their $T^\circ$-fixed vectors.  Such representations are interpreted as modules over an $\epsilon$-genuine $T^\circ$-spherical Hecke algebra (studied also in \cite[\S 3.2]{WWL1}).  Section \ref{Unramified} describes this Hecke algebra, in terms of the Brylinski-Deligne invariants functorially associated to $\alg{T}'$ in \cite{MWIntegral}.

Section \ref{Automorphic} returns to the global setting.  The commutativity of spherical Hecke algebras implies that irreducible admissible $\epsilon$-genuine representations of $\tilde T_\AA$ factor uniquely into $\epsilon$-genuine representations of $\tilde T_v$ at almost all places $v$ of a global field $F$.  Methods of Flath \cite{Flath} apply with minimal changes along the way.  In addition to studying generalities on the irreducible representations of $\tilde T_\AA$, this section defines the $\epsilon$-genuine automorphic representations of $\tilde T_\AA$ and classifies them in terms of their ``$\epsilon$-genuine automorphic central character.''  Section \ref{Automorphic} culminates in a multiplicity estimate for automorphic representations in $L_\chi^2(T_F \backslash \tilde T_\AA)$ for any cover of a torus over a global field.  For split tori, this is a \textbf{multiplicity-one} result, which should have applications to the global theory of Eisenstein series for covering groups.  A similar result, for some special cases in the function field setting, is found in very recent work of Lysenko \cite{Lysenko}.  The key is to estimate the index of $T_F Z(\tilde T_\AA)$ in a maximal abelian subgroup of $\tilde T_\AA$.  A result in harmonic analysis \cite{L-P} connects this index to the multiplicities in the automorphic spectrum.

Our results are heavily influenced by a desire to {\em parameterize} the irreducible admissible (or automorphic, in the global setting) representations of covering groups $\tilde G$ by something like Langlands parameters.  Here we believe that the L-group should be definable from a group $\alg{G}$, the Brylinski-Deligne invariants of an extension $\alg{K}_2 \Into \alg{G}' \Onto \alg{G}$, and a positive integer $n$ determining the degree of the cover.  Therefore, we have taken measures to describe everything possible in terms Brylinski-Deligne invariants as described in \cite{B-D} and \cite{MWIntegral}.  We plan to introduce such an L-group in a later article, and the work here should provide evidence for a more general conjecture.

\subsection*{Acknowledgments}

Our motivation for revisiting covers of tori came from a workshop ``Automorphic forms and harmonic analysis on covering groups'' at the American Institute of Mathematics, organized by Jeffrey Adams, Wee Teck Gan, and Gordan Savin.  We thank AIM for their support in this research, and many colleagues there for their insights.  In particular, Gordan Savin's work with Loke \cite{L-S} highlighted (for me) the importance of a multiplicity result for automorphic representations of covers of tori.  We thank Tamotsu Ikeda and Kaoru Hiraga for locating a crucial error in an early draft of this paper, and discussing other sections.  They were also generous to share some preliminary ideas related to endoscopy for covers of tori -- these provided some insights in studying the center of covers of tori over local fields.  Sergei Lysenko has worked on parallel results in the function field and geometric setting, and our conversations at AIM were most appreciated.

\section*{Preliminaries}

We work in this paper with smooth group schemes over fields, and occasionally over discrete valuation rings.  Our conventions follow \cite[Expos\'e VIII, IX]{SGA3} in our discussion of groups of multiplicative type (always of finite type over the base scheme).  In particular, when $F$ is a field, and $\alg{G}$ is a group of multiplicative type over $F$, then we write $\shom(\alg{G}, \alg{G}_{\mult})$ for its character lattice; this character lattice is naturally a local system of finitely-generated abelian groups on $F_{\et}$.  

An algebraic torus is a group of multiplicative type whose character lattice is a local system of finitely-generated free abelian groups.  When $\alg{T}$ is a torus, we write $\sheaf{Y}$ for its cocharacter lattice, and $\sheaf{X}$ for its character lattice, viewed as local systems on $F_{\et}$ in duality.

A \defined{quasitrivial} torus over $F$, sometimes called ``quasisplit'' or ``permutation'' or ``induced'' in the literature, will mean a torus over $F$ obtained by restriction of scalars from a split torus over a finite \'etale $F$-algebra.

When working with local systems on $F_{\et}$, we often write something like $y \in \sheaf{Y}$ to mean that $y \in \sheaf{Y}[E]$ for some finite \'etale $E/F$.  When describing constructions with \'etale sheaves, we often abuse notation in this way when our constructions are \'etale local.  

Our treatment of Galois cohomology mostly follows Serre, \cite{Serre}.  We frequently apply the \defined{Grunwald-Wang theorem} (discussed in, e.g., \cite[Chapter X]{ArtinTate}), which implies that $\Sha^1(\alg{\mu}_n) = 0$, whenever $F$ is a global field for which $\# \alg{\mu}_n(F) = n$.

Whenever $A$ is an abelian group, we write $A_{/n}$ for $A / n A$, and we write $A_{[n]}$ for the $n$-torsion subgroup of $A$.  Similar notation is applied for $n$-torsion in sheaves of abelian groups, group schemes, etc..

We work in the framework of Brylinski and Deligne \cite{B-D} when discussing central extensions of $\alg{T}$ by $\alg{K}_2$.  We refer to the first and second Brylinski-Deligne invariants, when referring to the results of \cite[Proposition 3.11]{B-D}, for example.  A description of these invariants, valid over a discrete valuation ring as well as a field, is found in our recent work \cite{MWIntegral}.

An \defined{lc-group} is a locally compact Hausdorff topological group.  An \defined{lca-group} is an abelian lc-group.  We work with continuous representations of lc-groups in this paper.  The cyclic group $\mu_n = \alg{\mu}_n(F)$ of order $n$ is used throughout the paper, and $\epsilon \From \mu_n \Into \CC^\times$ will be an injective character, fixed at all times.

If $\mu_n \Into \tilde G \Onto G$ is a central extension of lc-groups, then an \defined{$\epsilon$-genuine representation} of $\tilde G$ will mean a pair $(\pi, V)$ where $V$ is a complex vector space, $\pi \From \tilde G \To GL(V)$ is a homomorphism, and $\pi(\zeta) v = \epsilon(\zeta) \cdot v$ for all $\zeta \in \mu_n$.

We write $\Hom(\tilde G, \CC^\times)$ for the group of {\em continuous} homomorphisms.  We write $\Hom_\epsilon(\tilde G, \CC^\times)$ to be the subset of $\Hom(\tilde G, \CC^\times)$ consisting of characters which restrict to $\epsilon$ on $\mu_n$.  

\section{The center of a cover of a torus}
\label{Center}
\subsection{Tori and quadratic forms}

Let $\alg{T}$ be an algebraic torus over \textbf{any} field $F$.  Let $\alg{K}_2 \Into \alg{T}' \Onto \alg{T}$ be a central extension as studied by Brylinski and Deligne \cite{B-D}.  Let $n$ be a positive integer.

Write $\sheaf{X}$ and $\sheaf{Y}$ for the character and cocharacter lattices, viewed as local systems on $F_{\et}$.  These lattices are in duality, written $\langle \cdot, \cdot \rangle \From \sheaf{X} \otimes \sheaf{Y} \To \ZZ$.

Brylinski and Deligne associate to $\alg{T}'$ a quadratic form $Q \From \sheaf{Y} \To \ZZ$, i.e., an element $Q \in H_{\et}^0(F, \ssym^2 \sheaf{X})$.  Define a symmetric bilinear form $\beta \From \sheaf{Y} \otimes \sheaf{Y} \To n^{-1} \ZZ$ by
$$\beta(y_1, y_2) = n^{-1} \cdot \left( Q(y_1 + y_2) - Q(y_1) - Q(y_2) \right).$$
Define local systems $\sheaf{Y}^\sharp$ and $\sheaf{X}^\sharp$ by
$$n \sheaf{Y} \subset \sheaf{Y}^\sharp \defeq \{ y \in \sheaf{Y} : \beta(y, y') \in \ZZ \text{ for all } y' \in \sheaf{Y} \} \subset \sheaf{Y},$$
$$\sheaf{X} \subset \sheaf{X}^\sharp \defeq \{ x \in n^{-1} \sheaf{X} : \langle x, y \rangle \in \ZZ \text{ for all } y \in \sheaf{Y}^\sharp \} \subset n^{-1} \sheaf{X}.$$ 
The local systems $\sheaf{X}^\sharp$ and $\sheaf{Y}^\sharp$ are in duality via the pairing $\langle \cdot, \cdot \rangle$.

Let $\delta \From \sheaf{Y} \To n^{-1} \sheaf{X}$ denote the unique homomorphism of local systems satisfying
$$\langle \delta(y_1), y_2 \rangle = \beta(y_1, y_2) \text{ for all } y_1, y_2.$$
The duality between $\sheaf{X}$ and $\sheaf{Y}$, and between $\sheaf{X}^\sharp$ and $\sheaf{Y}^\sharp$, gives the following.
\begin{proposition}
\label{XQn}
$\delta(\sheaf{Y}) \subset \sheaf{X}^\sharp$ and $\delta(\sheaf{Y}^\sharp) \subset \sheaf{X}$.  In this way, we find an isomorphism,
$$\delta \From \frac{\sheaf{Y}}{\sheaf{Y}^\sharp} \xrightarrow{\sim} \frac{\sheaf{X}^\sharp}{\sheaf{X}}.$$
\end{proposition}
\proof
First, If $y \in \sheaf{Y}$ and $y^\sharp \in \sheaf{Y}^\sharp$, then $\langle \delta(y), y^\sharp \rangle = \beta(y, y^\sharp) \in \ZZ$.  Hence $\delta(y) \in \sheaf{X}^\sharp$.

Second, if $y^\sharp \in \sheaf{Y}^\sharp$ and $y \in \sheaf{Y}$, then $\langle \delta(y^\sharp), y \rangle = \beta(y^\sharp, y) \in \ZZ$.  Hence $\delta(y^\sharp) \in \sheaf{X}$.  

Finally, if $y \in \sheaf{Y}$ and $\delta(y) \in \sheaf{X}$, and $y' \in \sheaf{Y}$, then $\beta(y, y') = \langle \delta(y), y' \rangle \in \ZZ$.  Hence $y \in \sheaf{Y}^\sharp$.

These three observations imply that $\delta$ defines an injective homomorphism $\sheaf{Y} / \sheaf{Y}^\sharp \Into \sheaf{X}^\sharp / \sheaf{X}$.  Surjectivity follows from the fact that these are local systems of finite abelian groups in duality, hence locally of the same cardinality.  
\qed

\subsection{Related group schemes}

Define groups of multiplicative type over $F_{\et}$ according to their character lattices in the table below.
\begin{center}
\begin{tabular}{l|c|r}
Group $\alg{G}$ & Character lattice $\shom(\alg{G}, \alg{G}_{\mult})$ & Description \\ \hline
$\alg{T}$ & $\sheaf{X} \isom n^{-1} \sheaf{X}$ & Torus \\
$\alg{T}_{[n]}$ & $n^{-1} \sheaf{X} / \sheaf{X}$ & Finite $n$-torsion \\ 
$\alg{\hat T}$ & $\sheaf{Y} \isom n \sheaf{Y}$ & Torus \\
$\alg{\hat T}_{[n]}$ & $\sheaf{Y} / n \sheaf{Y}$ & Finite $n$-torsion \\ 
$\alg{T}^\sharp$ & $\sheaf{X}^\sharp$ & Torus \\
$\alg{\mu}$ & $\sheaf{X}^\sharp / \sheaf{X}$ & Finite $n$-torsion \\
$\alg{\hat \mu}$ & $\sheaf{Y} / \sheaf{Y}^\sharp$ & Finite $n$-torsion \\
$\alg{\nu}$ & $n^{-1} \sheaf{X} / \sheaf{X}^\sharp$  & Finite $n$-torsion \\ 
$\alg{\hat \nu}$ & $\sheaf{Y}^\sharp / n \sheaf{Y}$ & Finite $n$-torsion \\
\hline
\end{tabular}
\end{center}

Various maps of character lattices correspond contravariantly to $F$-homomorphisms of groups of multiplicative type.
\begin{center}
\begin{tabular}{c|c}
Group homomorphisms & Homomorphisms of character lattices  \\ \hline
$\alg{\mu} \Into \alg{T}^\sharp \xtwoheadrightarrow{i} \alg{T}$ & $\sheaf{X} \Into \sheaf{X}^\sharp \Onto \sheaf{X}^\sharp / \sheaf{X}$  \\
$\alg{\nu} \Into \alg{T} \xtwoheadrightarrow{j} \alg{T}^\sharp$ & $\sheaf{X}^\sharp \Into n^{-1} \sheaf{X} \Onto n^{-1} \sheaf{X} / \sheaf{X}^\sharp$  \\
$\alg{\nu} \Into \alg{T}_{[n]} \Onto \alg{\mu}$ & $\sheaf{X}^\sharp / \sheaf{X} \Into n^{-1} \sheaf{X} / \sheaf{X} \Onto n^{-1} \sheaf{X} / \sheaf{X}^\sharp$  \\ 
$\alg{\hat \mu} \Into \alg{\hat T}_{[n]} \Onto \alg{\hat \nu}$ & $\sheaf{Y}^\sharp / n \sheaf{Y} \Into \sheaf{Y} / n \sheaf{Y} \Onto \sheaf{Y} / \sheaf{Y}^\sharp$  \\ 
$\delta \From \alg{T}^\sharp \xrightarrow{} \alg{\hat T}$ & $\sheaf{Y} \xrightarrow{\delta} \sheaf{X}^\sharp$ \\
$\delta \From \alg{\mu} \xrightarrow{\sim} \alg{\hat \mu}$ & $\sheaf{Y} / \sheaf{Y}^\sharp \xrightarrow{\delta} \sheaf{X}^\sharp / \sheaf{X}$. \\
\hline
\end{tabular}
\end{center}

For what follows, it will be important to consider the following group of multiplicative type over $F$, and its associated character lattice:  define
\begin{equation}
\label{DefOfR}
\alg{R} = \{ (t, \hat t) \in \alg{T} \times \alg{\hat T} : \delta j(t) = \hat t^n \}.
\end{equation}
Write $p \From \alg{R} \To \alg{T}$ and $\hat p \From \alg{R} \To \alg{\hat T}$ for the maps given by $p(t, \hat t) = t$ and $\hat p(t, \hat t) = \hat t$.  The character lattice of $\alg{R}$ will be denoted $\sheaf{V}$,
$$\sheaf{V} = \frac{ \sheaf{X} \oplus \sheaf{Y} }{ \langle (n \delta(y), -ny) : y \in \sheaf{Y} \rangle }.$$
The map $p \From \alg{R} \To \alg{T}$ is surjective, since the corresponding map $x \mapsto (x,0)$ from $\sheaf{X}$ to $\sheaf{V}$ is injective.  It follows quickly that the free rank of $\sheaf{V}$ is equal to the rank of $\alg{T}$ (which equals the rank of $\sheaf{X}$, $\sheaf{Y}$, $\sheaf{X}^\sharp$, $\sheaf{Y}^\sharp$, etc.).  

Consider the homomorphism $e \From \sheaf{V} \To \sheaf{X}^\sharp$, $e(x,y) = \delta(y) + x$.  Proposition \ref{XQn} implies that $e$ is surjective, and thus defines an injective homomorphism $e \From \alg{T}^\sharp \Into \alg{R}$.  By consideration of rank, this identifies $\alg{T}^\sharp$ with the neutral component of $\alg{R}$.

The component group $\pi_0 \alg{R}$ is a finite group of multiplicative type, whose character lattice is identified with the kernel $\Ker(\sheaf{V} \xrightarrow{e} \sheaf{X}^\sharp)$.  Observe that, if $(x,y) \in \sheaf{X} \oplus \sheaf{Y}$ and $e(x,y) = 0$, then $x = - \delta(y)$.  It follows that $\delta(y) \in \sheaf{X}$, and so by Proposition \ref{XQn} we find that $y \in \sheaf{Y}^\sharp$.  We find that
$$\Ker(\sheaf{V} \xrightarrow{e} \sheaf{X}^\sharp) = \frac{ \{ (\delta(y), -y) : y \in \sheaf{Y}^\sharp \}}{\{ (n \delta(y), -ny) : y \in \sheaf{Y} \} } \isom \sheaf{Y}^\sharp / n \sheaf{Y}.$$
This identifies $\pi_0 \alg{R}$ with $\alg{\hat \nu}$, giving a short exact sequence
\begin{equation}
\label{SequenceTRv}
\alg{T}^\sharp \xhookrightarrow{e} \alg{R} \Onto \alg{\hat \nu}.
\end{equation}

The kernel of $i \From \alg{T}^\sharp \To \alg{T}$ is identified with $\alg{\mu}$.  Similarly, we compute the kernel of $p \From \alg{R} \To \alg{T}$.  The finite group scheme $\Ker(p)$ has character lattice
$$\frac{\sheaf{V} }{ \{ (x,0) : x \in \sheaf{X} \}} = \frac{\sheaf{X} \oplus \sheaf{Y}}{\sheaf{X} \oplus n\sheaf{Y}} \isom \frac{\sheaf{Y}}{n\sheaf{Y}}.$$
Thus we find a short exact sequence
\begin{equation}
\label{SequenceTRT}
\alg{\hat T}_{[n]} \Into \alg{R} \Onto \alg{T}.
\end{equation}
Chasing the character lattices around, the exact sequences \eqref{SequenceTRv} and \eqref{SequenceTRT} fit into a commutative diagram with exact rows and columns.
\begin{equation}
\label{RDiagram}
\begin{tikzcd}
\alg{\mu} \inarrow{r} \inarrow{d}{e}  & \alg{T}^\sharp \onarrow{r}{i} \inarrow{d}{e} & \alg{T} \arrow{d}{=} \\
\alg{\hat T}_{[n]} \inarrow{r} \onarrow{d} & \alg{R} \onarrow{r} \onarrow{d} & \alg{T} \\
\alg{\hat \nu} \arrow{r}{=} & \alg{\hat \nu}
\end{tikzcd}
\end{equation}

\subsection{Center for local fields}
\label{CenterCoverLocal}
Now assume that $F$ is a local field, and write $T = \alg{T}(F)$, $R = \alg{R}(F)$, etc..  Assume that $\# \mu_n(F) = n$.  Taking $F$-points of the extension $\alg{K}_2 \Into \alg{T}' \Onto \alg{T}$, and pushing out via the Hilbert symbol $\Hilb_n$ yields a central extension of locally compact groups (see \cite[Construction 10.3]{B-D})
\begin{equation}
\label{MetaTorus}
\mu_n \Into \tilde T \Onto T.
\end{equation}

The Kummer sequences for the tori $\alg{T}$, $\alg{T}^\sharp$, and $\alg{\hat T}$, yield a commutative diagram with exact rows.
\begin{equation}
\label{LocalKummer}
\begin{tikzcd}
T_{/n} \inarrow{r}{\Kum} \arrow{d}{j} & H_{\et}^1(F, \alg{T}_{[n]} ) \onarrow{r} \arrow{d}{j} & H_{\et}^1(F, \alg{T})_{[n]} \arrow{d}{j} \\
T_{/n}^{\sharp} \inarrow{r}{\Kum}  \arrow{d}{\delta} & H_{\et}^1(F, \alg{T}_{[n]}^\sharp ) \onarrow{r} \arrow{d}{\delta} & H_{\et}^1(F, \alg{T}^\sharp)_{[n]} \arrow{d}{\delta} \\
\hat T_{/n} \inarrow{r}{\Kum} & H_{\et}^1(F, \alg{\hat T}_{[n]} ) \onarrow{r} & H_{\et}^1(F, \alg{\hat T})_{[n]}
\end{tikzcd}
\end{equation}

As sheaves on $F_{\et}$, $\alg{T}_{[n]}$ represents $\sheaf{Y} \otimes \mu_n$ and $\alg{\hat T}_{[n]}$ represents $\sheaf{X} \otimes \mu_n$.  Tate local duality provides a perfect pairing:
\begin{equation}
\label{LocalTatePairing}
H_{\et}^1(F, \alg{T}_{[n]}) \times H_{\et}^1(F, \alg{\hat T}_{[n]}) \xrightarrow{\TD} \mu_n.
\end{equation}
Define a pairing by composing the Tate pairing with the Kummer coboundary,
$$\Hilb_T \From T_{/n} \times \hat T_{/n} \To \mu_n, \quad \Hilb_T(t, \hat t) = \TD( \Kum t, \Kum \hat t ).$$
We call this a Hilbert pairing, since for $\alg{T} = \alg{G}_{\mult}$ it coincides with the Hilbert symbol.
\begin{lemma}
If $H_{\et}^1(F, \alg{T})_{[n]} = H_{\et}^1(F, \alg{\hat T})_{[n]} = 0$, then the Hilbert pairing $\Hilb_T$ is nondegenerate.  In particular, if $\alg{T}$ is a quasitrivial torus, then $\Hilb_T$ is nondegenerate.
\end{lemma}
\proof
The conditions of the Lemma imply that the Kummer coboundaries $\Kum$ are isomorphisms, and so the Lemma follows directly from nondegeneracy of the Tate pairing.
\qed

From \cite[\S 7.5, \S 7.6, \S 9.8.3]{B-D}, the commutator of the extension \eqref{MetaTorus} is an alternating bilinear form $\Comm \From T \times T \To \mu_n$, given by 
\begin{equation}
\label{CommFormula}
\Comm(t_1,t_2) = \Hilb_T(\delta j  t_1,  t_2) = \Hilb_T( t_1, \delta j t_2)^{-1}.
\end{equation}
Define
$$Z^\dag(T) = \{ t_1 \in T : \Comm(t_1, t_2) = 1 \text{ for all } t_2 \in T \}.$$
The center of $\tilde T$ fits into a short exact sequence of abelian groups
$$\mu_n \Into Z(\tilde T) \Onto Z^\dag(T).$$

Tate local duality and \eqref{CommFormula} enable a description of $Z^\dag(T)$.
\begin{thm}
\label{CenterLocal}
There are containments:
$$\Im \left( T^\sharp \xrightarrow{i} T \right) \subset \Im \left( R \To T \right) \subset Z^\dag(T).$$
Each inclusion has finite index, bounded as below.
\begin{align*}
\# \left( \frac{ Z^\dag(T)}{ \Im ( R \To T ) } \right) & \leq \# H_{\et}^1(F, \alg{T})_{[n]}, \\
\# \left( \frac{ \Im ( R \To T )}{\Im ( T^\sharp \To T) } \right) & = \frac{ \# (R/T^\sharp) \cdot \# \mu}{\# \hat T_{[n]}} \leq \frac{ \# \nu \cdot \# \mu}{\# \hat T_{[n]}}.
\end{align*}
In particular, $Z^\dag(T) = \Im(R \xrightarrow{\iota} T)$ if $\alg{T}$ is an quasitrivial torus.  These containments are equalities if $\alg{T}$ is split.
\end{thm}
\proof
We begin by studying the inclusion
$$\Im \left( R \To T \right) \subset Z^\dag(T).$$
The commutator formula (\ref{CommFormula}) implies that
\begin{align*}
Z^\dag(T) &= \left\{ t_1 \in T : \Hilb_T(\delta j  t_1,  t_2) = 1 \text{ for all } t_2 \in T \right\}, \\
&= \left\{ t_1 \in T : \TD(\Kum \delta j  t_1,  \Kum t_2) = 1 \text{ for all } t_2 \in T \right\}.
\end{align*}
The Kummer sequences \eqref{LocalKummer} and Tate local duality give
\begin{align*}
\Im (R \To T) &= \{ t_1 \in T : \delta j t_1 = \hat t^n \text{ for some } \hat t \in \hat T \}, \\
&= \{ t_1 \in T : \Kum \delta j t_1 = 1 \}, \\
&= \{ t_1 \in T : \TD(\Kum \delta j t_1 , \eta) = 1 \text{ for all } \eta \in H_{\et}^1(F, \alg{T}_{[n]}) \}. \\
\end{align*}

The map $t_1 \mapsto \TD(\Kum \delta j t_1, \bullet)$ and the Kummer sequence for $\alg{T}$ give an injective homomorphism
$$\frac{ Z^\dag(T)}{ \Im ( R \To T )} \Into \Hom \left( \frac{H_{\et}^1(F, \alg{T}_{[n]})}{\Kum(T_{/n})}, \mu_n \right) \isom \Hom \left( H_{\et}^1(F, \alg{T})_{[n]}, \mu_n \right).$$
This demonstrates the first part of the theorem.  Also, if $\alg{T}$ is quasitrivial then $H_{\et}^1(F, \alg{T}) = 0$ and so $Z^\dag(T) = \Im(R \To T)$ as claimed.

Diagram \eqref{RDiagram} identifies
$$\frac{ \Im ( R \To T )}{\Im ( T^\sharp \To T) } \isom \frac{R / \hat T_{[n]}}{T^\sharp / \mu} \isom \frac{ R/T^\sharp}{\hat T_{[n]} / \mu}.$$ 
As $R / T^\sharp$ is identified with a subgroup of $\nu$, the second part of the Theorem follows.

When $\alg{T}$ is a split torus, the sequence $\mu \Into \hat T_{[n]} \Onto \nu$ is exact, by considerations of cardinality alone.  It follows that, when $\alg{T}$ is split, we have
$$\Im(R \To T) = \Im(T^\sharp \To T),$$
finishing the proof of the Theorem.
\qed

\begin{remark}
Our earlier paper \cite{MWTori} demonstrated that $Z^\dag(T) = \Im(T^\sharp \To T)$ when $\alg{T}$ is split.
\end{remark}

The isogeny (on $F$-points) $T^\sharp \xrightarrow{i} T$ plays a critical role in our parameterization of representations of $\tilde T$.  For this reason, we give its image a name.
\begin{definition}
Define $C^\dag(T) = \Im(T^\sharp \To T)$.  The \defined{central core} of $\tilde T$ is defined to be the preimage $C(\tilde T)$ of $C^\dag(T)$ in $\tilde T$:
$$\mu_n \Into C(\tilde T) \Onto C^\dag(T).$$
\end{definition}
The central core of $\tilde T$ lies inside the center of $\tilde T$.  The previous Theorem estimates the index.
$$\# \left( \frac{Z(\tilde T)}{C(\tilde T)} \right) \leq  \# H_{\et}^1(F, \alg{T})_{[n]} \cdot \left( \frac{ \# \nu \cdot \# \mu}{\# \hat T_{[n]}} \right).$$

\subsection{Center for global fields}

Now suppose that $F$ is a global field, $\alg{T}$ is a torus over $F$, and $\alg{K}_2 \Into \alg{T}' \Onto \alg{T}$ is a central extension.  Suppose that $\# \mu_n(F) = n$.  The construction of \cite[\S 10.4]{B-D} gives a central extension
$$\mu_n \Into \tilde T_\AA \Onto T_\AA,$$
canonically split over $T_F$.

Let $\VV$ be the set of places of $F$, and for every place $v \in \VV$, write $T_v = \alg{T}(F_v)$.  Let $\mu_n \Into \tilde T_v \Onto T_v$ be the central extension obtained from $\tilde T_\AA$ by pulling back, and let $\tilde Z_v$ denote the center of $\tilde T_v$.  The extension $\tilde T_v$ is naturally identified with the the extension obtained by \cite[Construction 10.3]{B-D} and the Hilbert symbol, as in the last section.  At almost all places $v \in \VV$ (excluding all archimedean places and those where a splitting field of $\alg{T}$ ramifies, among others), the construction of \cite[\S 10.7]{B-D} provides a splitting $T_v^\circ \Into \tilde T_v$ over the maximal compact subgroup $T_v^\circ \subset T_v$.  

Let $\tilde Z_\AA$ denote the center of $\tilde T_\AA$.  The center of $\tilde T_v$ fits into a short exact sequence
$$\mu_n \Into \tilde Z_v \Onto Z_v^\dag,$$
where we write $Z_v^\dag$ to abbreviate $Z^\dag(T_v)$.  The commutator pairing on $T_\AA$ is the product of local commutator pairings:
$$\Comm(x,y) = \prod_{v \in \VV} \Comm(x_v, y_v),$$
for all $x,y \in T_\AA$.  Hence the center of $\tilde T_\AA$ fits into a short exact sequence
$$\mu_n \Into \tilde Z_\AA \Onto Z_\AA^\dag,$$
where
\begin{equation}
\label{LocalDescriptionCenter}
Z_\AA^\dag = \{ (t_v)_{v \in \VV} \in T_\AA : t_v \in Z_v^\dag \text{ for all } v \in \VV \}.
\end{equation}
This identifies $Z_\AA^\dag$ as the restricted direct product of the groups $Z_v^\dag$ with respect to their open compact subgroups $Z_v^\circ = Z_v^\dag \cap T_v^\circ$.

The center of $\tilde T_\AA$ can be characterized analogously to Theorem \ref{CenterLocal}.
\begin{thm}
\label{GlobalCenterImage}
There are containments:
$$\Im \left( T_\AA^\sharp \To T_\AA \right) \subset \Im \left( R_\AA \To T_\AA \right) \subset Z_\AA^\dag.$$
If $H_{\et}^1(F_v, \alg{T})_{[n]} = 0$ for all places $v$ (in particular, if $\alg{T}$ is a quasitrivial $F$-torus), then
$$Z_\AA^\dag = \Im \left( R_\AA \rightarrow T_\AA \right).$$
If $\alg{T}$ is split over $F$, then
$$Z_\AA^\dag = \Im \left( T_\AA^\sharp \To T_\AA \right).$$
\end{thm}
\proof
By the local description of $Z_\AA^\dag$ in \eqref{LocalDescriptionCenter}, and the local containments of Theorem \ref{CenterLocal}, we find the local containments of the Theorem:
$$\Im \left( T_\AA^\sharp \To T_\AA \right) \subset \Im \left( R_\AA \rightarrow T_\AA \right) \subset Z_\AA^\dag.$$

For what follows, choose models of $\alg{R}$ and $\alg{T}$ and the homomorphism $\alg{R} \To \alg{T}$, over the $S$-integers $\OO_S \subset F$, for some finite set of places $S$.  For $v \in \VV - S$, we obtain compact open subgroups $R_v^\circ = \alg{R}(\OO_v) \subset R_v = \alg{R}(F_v)$, and $T_v^\circ = \alg{T}(\OO_v) \subset T_v = \alg{T}(F_v)$.  Choosing different integral models of $\alg{R}$ and $\alg{T}$ will only change $R_v^\circ$ and $T_v^\circ$ at a finite set of places.

Suppose that $H_{\et}^1(F_v, \alg{T})_{[n]} = 0$ for all places $v \in \VV$, and $t = (t_v) \in Z_\AA^\dag$.  For all $v \in \VV$, there exists (by Theorem \ref{CenterLocal}) an element $r_v \in R_v$ which maps to $t_v$ via $\alg{R} \To \alg{T}$.  It remains to check that the resulting family of elements $(r_v)$ is an element of $R_\AA$.  In other words, we must check that $r_v \in R_v^\circ$ for almost all places $v \in \VV$.  

The algebraic homomorphism $p \From \alg{R} \To \alg{T}$ corresponds to an injective map of character groups $p \From \sheaf{X} \hookrightarrow \sheaf{V}$.  Let $E$ be a finite Galois extension of $F$, such that $\sheaf{X}$ and $\sheaf{V}$ restrict to constant sheaves on $E_{\et}$.  Let $w$ be a place of $E$ lying above $v$, with $v \not \in S$ nonarchimedean place, and write $\OO_w$ for the ring of integers in $E_w$.  

\'Etale descent characterizes the subgroups $T_v^\circ \subset T_v$ and $R_v^\circ \subset R_v$:
$$R_v^\circ =  \{ s \in R_v : \phi(s) \in \OO_w^\times \text{ for all } \phi \in \sheaf{V}[E] \};$$
$$T_v^\circ = \{ t \in T_v : \xi(s) \in \OO_w^\times \text{ for all } \xi \in \sheaf{X}[E] \}.$$

If $r_v \in R_v$ maps to $t_v \in T_v^\circ$, then $\xi(t_v) = p(\xi)(r_v) \in \OO_w^\times$ for all $\xi \in \sheaf{X}[E]$.  But the cokernel of $p \From \sheaf{X}[E] \Into \sheaf{V}[E]$ is finite, so for all $\phi \in \sheaf{V}[E]$, $\phi(r_v)^N \in \OO_w^\times$ for some positive integer $N$.  Since $E_w^\times / \OO_w^\times \isom \ZZ$, this implies that $\phi(r_v) \in \OO_w^\times$ for all $\phi \in \sheaf{V}[E]$.  Hence $r_v \in R_v^\circ$.

Since $t_v \in T_v^\circ$ for almost all places, we find that $r_v \in R_v^\circ$ at almost all places.  Thus $r = (r_v) \in R_\AA$.  

In the split case, the same arguments apply, replacing $\alg{R}$ by $\alg{T}^\sharp$ as needed.
\qed

\section{Irreducible representations, local case}
\label{Irreps}
Let $\alg{T}$ be a torus over a local field $F$ and consider a covering 
$$\mu_n \Into \tilde T \Onto T$$
as in Section \ref{CenterCoverLocal}.  There is a unique maximal compact subgroup $T^\circ \subset T$, and its preimage $\tilde T^\circ$ is a maximal compact subgroup of $\tilde T$.

The commutator pairing descends to a nondegenerate alternating bilinear form (on a finite $n$-torsion abelian group)
$$\Comm: \frac{T}{Z^\dag(T)} \times \frac{T}{Z^\dag(T)} \rightarrow \mu_n.$$
It follows that there is a \defined{Lagrangian decomposition} $L \times L^\ast =T / Z^\dag(T)$; this means that $L, L^\ast$ are subgroups of $T / Z^\dag(T)$ such that $\Comm \From L \times L^\ast \To \mu_n$ identifies $L^\ast$ with the Pontrjagin dual of $L$, $L \cdot L^\ast = T / Z^\dag(T)$, and 
\begin{align*}
\Comm( x,y) &= 0 \text{ for all } x,y \in L, \\
\Comm(x,y) &= 0 \text{ for all } x,y \in L^\ast.\\
\end{align*}

In particular, the quotient $T / Z^\dag(T)$ has order $(\# L)^2$.  The order of such a Lagrangian subgroup $L$ will be called the \defined{central index} of $\tilde T$, written:
$$\zind(\tilde T) \defeq \sqrt{ \# (T /  Z^\dag(T))} = \# L.$$

When $K$ is a compact lc-group, an \defined{irrep} of $K$ will mean an irreducible continuous representation of $K$ on a finite-dimensional complex vector space.  We write $\Irr(K)$ for the set of equivalence classes of such irreps.

If $K \To GL(V)$ is a representation of $K$ on a complex vector space, and $[\rho,W] \in \Irr(K)$ is an equivalence class of irreps, then we write $V_{[\rho]}$ for the isotypic subspace,
$$V_{[\rho]} = \Span_\CC \{ \Im(f) : f \in \Hom_K(W, V) \}.$$

In what follows, we define {\em irreducible admissible $\epsilon$-genuine} representations of $\tilde T$, construct them by means of Lagrangian subgroups, and classify them by central character.  Much of this is ``well-known'' and we mention \cite{MWTori} and \cite{ABPTV} and \cite{K-P} as other locations where such results appear.

\subsection{Admissible representations}

We define admissible representations of $\tilde T$, beginning in the nonarchimedean case.  Recall that $\tilde T^\circ$ is a maximal compact subgroup of $\tilde T$.
\subsubsection{Nonarchimedean case}
\begin{definition}
Suppose that $F$ is nonarchimedean, so that $\tilde T$ is a totally disconnected lc-group.  An \defined{admissible representation} of $\tilde T$ means a representation $\pi \From \tilde T \To GL(V)$ such that $V = \bigoplus_\rho V_{[\rho]}$ is the direct sum of its $\tilde T^\circ$-isotypic subspaces, and each $V_{[\rho]}$ is finite-dimensional.  An admissible representation $(\pi, V)$ is called \defined{irreducible} if its only $\tilde T$-stable subspaces are $\{ 0 \}$ and $V$ itself.
\end{definition} 

\begin{remark}
In the nonarchimedean setting, every admissible representation $(\pi, V)$ of $\tilde T$ is \defined{smooth}:  every vector $v \in V$ has open stabilizer in $\tilde T$.  Every {\em irreducible} smooth representation of $\tilde T$ is admissible.  We use admissible representations mostly for consistency with what comes later.
\end{remark}

\subsubsection{Archimedean case}
Now suppose that $F$ is archimedean, so that $\tilde T \To T$ is a covering of Lie groups.  As a topological group, $T$ is isomorphic to a product of some copies of $\RR^\times$, some copies of $\CC^\times$, and some copies of $U(1)$.  There are a few different ``brands'' of representation theory for noncompact Lie groups, but we demonstrate here that the ``Harish-Chandra module'' approach coincides with the straightforward ``continuous finite-dimensional approach'' when it comes to representations of $\tilde T$.

Write $N$ for the neutral component of $T$ and $\tilde N$ for the neutral component of $\tilde T$.  \textbf{Beware} that $\tilde N$ may not coincide with the preimage of $N$ in $\tilde T$, although the projection $\tilde T \To T$ restricts to a surjective map $\tilde N \To N$.  Write $N^\circ = N \cap T^\circ$ and $\tilde N^\circ = \tilde N \cap \tilde T^\circ$.  These are maximal compact subgroups of $N$ and $\tilde N$, respectively. 

Note that $\tilde N^\circ$ and $N^\circ$ are honest tori -- products of circles -- and $\tilde N^\circ \To N^\circ$ is a cover of finite degree (dividing $n$).    By connectedness and the continuity of the commutator, we find that $\tilde N$ is contained in the center of $\tilde T$.  Since $T = T^\circ \cdot N$, we have $\tilde T = \tilde T^\circ \cdot \tilde N$ as well.  More precisely, multiplication gives an isomorphism
$$\tilde T^\circ \times_{\tilde N^\circ} \tilde N \defeq \frac{ \tilde T^\circ \times \tilde N}{ \langle (\tilde u, 1) = (1, \tilde u) : \tilde u \in \tilde N^\circ \rangle } \xrightarrow{\sim} \tilde T.$$

Define $\Lie{t}$ to be the complexified Lie algebra of $T$, and $\Lie{t}^\circ$ to be the complexified Lie algebra of $T^\circ$.  The cover $\tilde N \To N$ identifies $\Lie{t}$ with the complexified Lie algebra of $\tilde N$ (hence of $\tilde T$), and identifies $\Lie{t}^\circ$ with the complexified Lie algebra of $\tilde N^\circ$.

\begin{definition}
An \defined{admissible representation} of $\tilde T$ means a $(\Lie{t}, \tilde T^\circ)$-module $V$, such that $V = \bigoplus_{\rho} V_{[\rho]}$ is the direct sum of its $\tilde T^\circ$-isotypic subspaces, and each $V_{[\rho]}$ is finite-dimensional.  An admissible representation $V$ of $\tilde T$ is called \defined{irreducible} if its only $(\Lie{t}, \tilde T^\circ)$-stable subspaces are $\{ 0 \}$ and $V$ itself.
\end{definition}
\begin{remark}
Here, an $(\Lie{t}, \tilde T^\circ)$-module is a special case of the general notion of a $(\Lie{g}, K)$-module.  In other words, it requires simultaneous actions $\Lie{t} \To \End(V)$ and $\tilde T^\circ \To GL(V)$, such that the differential of the latter coincides with the former after restriction to $\tilde T^\circ$.  The condition on $(\Lie{g}, K)$-modules that reads $(Ad(k) X) v = k X k^{-1} v$ (for $X \in \Lie{g}$, $k \in K$) corresponds here to the condition that the $\Lie{t}$-action and $\tilde T^\circ$-action commute.
\end{remark}

Consider now a finite-dimensional continuous irreducible representation $\pi \From \tilde T \To GL(V)$.  Since $\tilde N \subset Z(\tilde T)$ and $\tilde T = \tilde T^\circ \cdot \tilde N$, we find that $(\pi, V)$ is irreducible when restricted to its maximal compact subgroup $\tilde T^\circ$.  Moreover, we find by Schur's lemma that $\tilde N$ acts via a continuous (hence smooth) character $\chi \From \tilde N \To \CC^\times$.  This yields the structure of a $(\Lie{t}, \tilde T^\circ)$-module on $V$, with $\Lie{t}$ acting by $d \chi$.

Conversely, given an irreducible $(\Lie{t}, \tilde T^\circ)$-module $V$, the commutativity of $\Lie{t}$ (and the fact that its action commutes with $\tilde T^\circ$) implies that $V$ is irreducible as a representation of $\tilde T^\circ$ alone.  The coincidence of the $\Lie{t}$ action and the $\tilde T^\circ$ action after restriction to the subgroup $\tilde N^\circ$ implies that the $\Lie{t}$ action is the derivative of a character of $\tilde N$.  Putting this together, we find the following result.  
\begin{proposition}
The construction above gives a bijection between the sets of...
\begin{itemize}
\item
equivalence classes of finite-dimensional continuous irreducible representations $\pi \From \tilde T \To GL(V)$;
\item
equivalence classes of irreducible admissible representations of $\tilde T$.
\end{itemize}
\end{proposition}

Hereafter, an irreducible admissible representation of $\tilde T$ (by definition, an $(\Lie{t}, \tilde T^\circ)$-module) will be thought of as a finite-dimensional continuous irreducible representation of the group $\tilde T$.

\subsubsection{General classification}
Now take $F$ to be a local field, archimedean or nonarchimedean.  An admissible representation of $\tilde T$ is called unitary, if there exists a positive-definite $\tilde T$-invariant Hermitian form on the space of the representation.

\begin{definition}
Define $\Irr_\epsilon(\tilde T)$ to be the set of equivalence classes of irreducible $\epsilon$-genuine admissible representations of $\tilde T$.  Define $\Irr_\epsilon^{\unit}(\tilde T)$ to be the subset of equivalence classes of unitary representations.
\end{definition}

A description of $\Irr_\epsilon(\tilde T)$ in the nonarchimedean case follows from \cite[Theorem 3.1]{MWTori} (based on previous work of many others).  This is an analogue of the far older Stone von-Neumann theorem.  The archimedean case can be found in \cite[Proposition 2.2]{ABPTV}.
\begin{thm}
\label{LocalSvN}
Let $\chi \From Z(\tilde T) \To \CC^\times$ be an $\epsilon$-genuine continuous character.  Choose a Lagrangian decomposition $L \times L^\ast = T / Z^\dag(T) = \tilde T / Z(\tilde T)$, and let $\tilde M$ be the preimage of $L$ in $\tilde T$ (so $\tilde M$ is a maximal commutative subgroup of $\tilde T$).  Let $\chi_M$ be an extension of $\chi$ to $\tilde M$ and let $\pi = \Ind_{\tilde M}^{\tilde T} \chi_M$ be the (algebraically) induced representation.

Then $\pi$ is an irreducible admissible $\epsilon$-genuine representation of $\tilde T$ of dimension $\zind(\tilde T)$ and central character $\chi$, and its equivalence class does not depend on the choice of Lagrangian decomposition nor on the choice of extension of $\chi$.  If $\pi'$ is any irreducible representation of $\tilde T$ of central character $\chi$, then $\pi'$ is equivalent to $\pi$.  The representation $\pi$ is unitary if and only if the character $\chi$ is unitary.
\end{thm}

\begin{remark}
Strictly speaking, \cite[Proposition 2.2]{ABPTV} does not mention the construction of irreducible representations by induction from a polarizing subgroup.  But this follows easily from examination of the character below.
\end{remark}

This theorem determines a bijection 
$$\Irr_\epsilon(\tilde T) \leftrightarrow \Hom_\epsilon(Z(\tilde T), \CC^\times).$$
For $\chi \in \Hom_\epsilon(Z(\tilde T), \CC^\times)$, write $[\pi_\chi]$ for the equivalence class of $\Ind_{\tilde M}^{\tilde T} \chi_M$ defined above.

The representation $[\pi_\chi]$ is a finite-dimensional representation, induced from a character of a finite-index subgroup.  As such, its character is straightforward to compute.
\begin{proposition}
\label{CharacterIrrepT}
The character of $[\pi_\chi]$ is given by
$$\Tr [\pi_\chi](\tilde t) = \begin{cases} \zind(\tilde T) \cdot \chi(\tilde t) & \text{ if } \tilde t \in Z(\tilde T); \\ 0 & \text{ otherwise }. \end{cases}$$
\end{proposition}
\proof
If $\tilde t \in Z(\tilde T)$, then $\pi_\chi(\tilde t)$ is the scalar operator $\chi(\tilde t)$, operating on a space of dimension $\zind(\tilde T)$.  This verifies the proposition when $\tilde t \in Z(\tilde T)$.

When $\tilde t \not \in Z(\tilde T)$, it projects to an element $t \in T$ with $t \not \in Z^\dag(T)$.  There exists a Lagrangian decomposition $L \times L^\ast = T / Z^\dag(T)$, such that $L$ contains $t \text{ mod } Z^\dag(T)$.  Let $\tilde M$ be the preimage of $L$ in $\tilde T$.  The projection $\tilde M \Onto L$ and the identification $L^\ast = \Hom(L, \mu_n)$, allows every element $\lambda^\ast \in L^\ast$ to be viewed as a character of $\tilde M$.  

Fixing one extension $\chi_M \in \Hom(\tilde M, \CC^\times)$ of $\chi$, the restriction of $\pi_\chi$ to $\tilde M$ is a direct sum of characters:
$$\Res_{\tilde M} \pi_\chi = \bigoplus_{\lambda^\ast \in L^\ast} \chi_M \cdot \lambda^\ast.$$
Since $\tilde t$ projects to a nontrivial element of $L$,
$$\Tr [\pi_\chi](\tilde t) = \chi_M(\tilde t) \cdot \sum_{\lambda^\ast \in L^\ast} \lambda^\ast (\tilde t) = 0.$$
\qed 

\subsection{Central core character and pouches}

The central character gives a bijection,
$$\Irr_\epsilon(\tilde T) \leftrightarrow \Hom_\epsilon( Z(\tilde T), \CC^\times).$$
However, the center of $\tilde T$ is difficult to describe in general.  Instead, we have described the subgroup
$$C(\tilde T) = \Im \left( \tilde T^\sharp \To \tilde T \right) \subset Z(\tilde T),$$
called the central core.

\subsubsection{Central core character}
Since $C(\tilde T)$ is a subgroup of $Z(\tilde T)$ of finite index, restriction gives a surjective map
$$\Hom_\epsilon(Z(\tilde T), \CC^\times) \Onto \Hom_\epsilon(C(\tilde T), \CC^\times),$$
whose fibres have constant cardinality $\# \left(Z^\dag(T) / C^\dag(T) \right)$.  Define the \defined{core index} of $\tilde T$ by
$$\cind(\tilde T) = \# \left(Z^\dag(T) / C^\dag(T) \right).$$
\begin{definition}
Suppose that $[\pi] \in \Irr_\epsilon(\tilde T)$.  The \defined{central core character} of $[\pi]$ is defined to be the character of $C(\tilde T)$ obtained by taking the central character and restricting to $C(\tilde T)$.  This gives a finite-to-one surjective map
$$\Irr_\epsilon(\tilde T) \To \Hom_\epsilon(Z(\tilde T), \CC^\times) \To \Hom_\epsilon(C(\tilde T), \CC^\times).$$
We say that two equivalence classes $[\pi_1], [\pi_2] \in \Irr_\epsilon(\tilde T)$ belong to the same \defined{pouch} if they have the same central core character.
\end{definition}
In this way, $\Irr_\epsilon(\tilde T)$ is partitioned into pouches, and each pouch has cardinality equal to the core index $\cind(\tilde T)$.  For all $\xi \in \Hom_\epsilon(C(\tilde T), \CC^\times)$, write $\Pouch_\xi \subset \Irr_\epsilon(\tilde T)$ for the pouch with central core character $\xi$.

\begin{corollary}
For all $\xi \in \Hom_\epsilon(C(\tilde T), \CC^\times)$,
$$\sum_{[\pi] \in \Pouch_\xi} \Tr[\pi] (\tilde t) = \begin{cases} \cind(\tilde T) \cdot \zind(\tilde T)\cdot  \xi(\tilde t) & \text{ if } \tilde t \in C(\tilde T); \\
0 & \text{ otherwise.} \end{cases}
$$
\end{corollary}
\proof
This follows directly from Proposition \ref{CharacterIrrepT}, and standard facts about characters of finite abelian groups.
\qed

While we like to put irreducible representations together into pouches, for later purposes of parameterization, it is useful to parameterize the elements of the pouches so as not to forget the individual representations.  

Define the group $P(T) = Z^\dag(T) / C^\dag(T)$ with cardinality $\cind(\tilde T)$.  Pushing out via $\xi$ yields a group $P_\xi(T) = Z(\tilde T) / \Ker(\xi)$ fitting into commutative diagram with exact rows.
$$\begin{tikzcd}
C(\tilde T) \inarrow{r} \arrow{d}{\xi} & Z(\tilde T) \onarrow{r} \arrow{d} & P(T) \arrow{d}{=} \\
\CC^\times \inarrow{r} & P_\xi(T) \onarrow{r} & P(T).
\end{tikzcd}$$

The elements of $\Pouch_\xi$ are in bijection wth characters of $Z(\tilde T)$ extending $\xi$.  To give a character of $Z(\tilde T)$ extending $\xi$ is the same as giving a splitting of the bottom row of this commutative diagram.  Such splittings may be thought of as ``twisted characters'' of the finite group $P(T)$, and they parameterize the elements of the pouch $\Pouch_\xi$.  

\subsubsection{Pulling back to $\tilde T^\sharp$}

Recall that the central core $C(\tilde T)$ is the image of the map $\tilde T^\sharp \To \tilde T$.  Here, $\tilde T^\sharp$ fits into an \textbf{abelian} extension
$$\mu_n \Into \tilde T^\sharp \Onto T^\sharp,$$
which also arises from the construction of Brylinski-Deligne.  An $\epsilon$-genuine character $\xi \From C(\tilde T) \To \CC^\times$ pulls back to an $\epsilon$-genuine character $\phi \From \tilde T^\sharp \To \CC^\times$.  

Putting it all together, we find a map
$$\Irr_\epsilon(\tilde T) \Onto \Hom_\epsilon(C(\tilde T), \CC^\times) \Into \Hom_\epsilon(\tilde T^\sharp, \CC^\times) = \Irr_\epsilon(\tilde T^\sharp).$$

The kernel of $T^\sharp \xrightarrow{i} T$ is the finite group $\mu = \alg{\mu}(F)$.  We expect that in a Langlands-style parameterization of $\Irr_\epsilon(\tilde T)$, the following intermediate theorem will be useful.  It also serves to summarize the results of this section.
\begin{thm}
\label{ParamLocalTori}
Let $\alg{T}$ be a torus over a local field.  Let $\alg{K}_2 \Into \alg{T}' \Onto \alg{T}$ be a central extension, and $n \geq 1$ a positive integer for which $\# \mu_n(F) = n$.  Let $\mu_n \Into \tilde T \Onto T$ be the resulting central extension.  Let $i \From \alg{T}^\sharp \To \alg{T}$ the resulting isogeny with kernel $\alg{\mu}$.  The map $[\pi] \To \xi \circ i$, sending an element of $\Irr_\epsilon(\tilde T)$ to the pullback (via $i$) of its central core character $\xi$, defines a finite-to-one function
$$\Phi \From \Irr_\epsilon(\tilde T) \rightarrow \Hom_\epsilon(\tilde T^\sharp, \CC^\times) = \Irr_\epsilon(\tilde T^\sharp).$$
The image of this finite-to-one map consists of those $\epsilon$-genuine characters of $\tilde T^\sharp$ which are trivial on $\mu$.  The nonempty fibres, called pouches, have cardinality equal to $\cind(\tilde T)$, which is bounded by $\# H_{\et}^1(F, \alg{T})_{[n]} \cdot \# \nu \cdot \# \mu / \# \hat T_{[n]}$.  In particular, if $\alg{T}$ is split, then the pouches are singletons and $\Irr_\epsilon(\tilde T)$ embeds into $\Irr_\epsilon(\tilde T^\sharp)$.
\end{thm}

\section{Unramified representations}
\label{Unramified}
In this section, let $\alg{T}$ be a torus over $\OO$, the ring of integers in a nonarchimedean local field $F$ with residue field $\FF_q$.  Let $T^\circ = \alg{T}(\OO)$ and $T = \alg{T}(F)$.  The valuation on $F$ will always be normalized so that $\val \From F^\times / \OO^\times \xrightarrow{\sim} \ZZ$ is an isomorphism.  Let $n$ be a positive integer which divides $q-1$; thus $\mu_n(\FF_q)$ is a cyclic group of order $n$, identified via the Teichmuller lift with $\mu_n(F)$.  If $z \in \OO$, we write $\bar z$ for its projection in $\FF_q$.

It will be helpful to choose an unramified extension $E/F$ over which $\alg{T}$ splits.  Write $\OO_E$ for the ring of integers in $E$.  Let ${\Fr}$ be a generator of $\Gal(E/F)$.  The \'etale local systems $\sheaf{X}$ and $\sheaf{Y}$ are encoded in the $\ZZ[{\Fr}]$-modules $X = \sheaf{X}[E]$ and $Y = \sheaf{Y}[E]$, in ${\Fr}$-invariant duality
$$\langle \cdot, \cdot \rangle  \From X \times Y \To \ZZ, \quad \langle {\Fr}(x), {\Fr}(y) \rangle = \langle x,y \rangle.$$
As $\alg{T}$ splits over $E$, we may identify Galois-fixed characters and cocharacters with those fixed by Frobenius:
$$Y^{\Fr} = H_{\et}^0(F, \sheaf{Y}), \quad X^{\Fr} = H_{\et}^0(F, \sheaf{X}).$$
There is a unique isomorphism $\val_T \From T/T^\circ \xrightarrow{\sim} Y^{\Fr}$ such for all $\xi \in X$ and all $t \in T$,
$$\val(\xi(t)) = \langle \xi, \val_T(t) \rangle.$$ 

Throughout this section, let $\alg{K}_2 \Into \alg{T}' \Onto \alg{T}$ be a central extension defined over $\OO$.  Such central extensions have been described and classified in \cite{MWIntegral}.  This gives a central extension
$$\mu_n \Into \tilde T \Onto T,$$
together with a splitting $T^\circ \Into T$ over the maximal compact subgroup $T^\circ$ (see also \cite[\S 10.7]{B-D}).

Brylinski and Deligne associate two invariants to the extension $\alg{K}_2 \Into \alg{T}' \Onto \alg{T}$.  We describe these invariants (in the $\OO$-integral case) in \cite{MWIntegral}.  The first invariant gives a ${\Fr}$-invariant quadratic form $Q \From Y \To \ZZ$.  The second invariant is a central extension of sheaves of groups on $\OO_{\et}$,
$$\sheaf{G}_m \Into \sheaf{D} \Onto \sheaf{Y}.$$
Writing $D^\circ = \sheaf{D}[\OO_E]$, this extension gives a $\Fr$-equivariant central extension of groups,
$$\OO_E^\times \Into D^\circ \Onto Y.$$

\subsection{Genuine spherical representations}
The theory of $\epsilon$-genuine spherical representations of $\tilde T$ has been developed to some extent in \cite[\S 5, \S 6]{MWTori}, \cite[\S 3.2]{WWL1}, and \cite{Lysenko}.  We go further here, exploiting some of the structural results of \cite{MWIntegral}.  But first, we recall the basic definitions.

\begin{definition}
An $\epsilon$-genuine \defined{spherical} representation of $\tilde T$ is an $\epsilon$-genuine smooth representation which is generated by its $T^\circ$-fixed vectors.  With $\tilde T$-intertwining maps as morphisms, these form an abelian category $\Cat{Rep}_\epsilon(\tilde T; T^\circ)$.
\end{definition}

\begin{definition}
The $\epsilon$-genuine \defined{spherical Hecke algebra} is the set $\hecke_\epsilon(\tilde T; T^\circ)$ of compactly-supported functions $f \From \tilde T \To \CC$ satisfying
$$f(\zeta k_1 \tilde t k_2) = \epsilon(\zeta) f(\tilde t) \text{ for all } \zeta \in \mu_n, k_1, k_2 \in T^\circ, \tilde t \in \tilde T.$$
Fix the Haar measure on $T$ for which $T^\circ$ has measure $1$.  Convolution in the Hecke algebra is given by
$$[f_1 \ast f_2](t) = \int_{T} f_1(x) \cdot f_2(x^{-1} t) dx.$$
(The reader may check that the integrand is a well-defined function on $T = \tilde T / \mu_n$.)
\end{definition}

When $(\pi, V)$ is a $\epsilon$-genuine spherical representation of $\tilde T$, the set $V^{T^\circ}$ of $T^\circ$-fixed vectors forms a module over the Hecke algebra $\hecke_\epsilon(\tilde T; T^\circ)$.  This defines an equivalence of abelian categories:
$$\Cat{Rep}_\epsilon(\tilde T; T^\circ) \xrightarrow{\sim} \Cat{Mod} \left( \hecke_\epsilon(\tilde T; T^\circ) \right).$$

\subsection{Spherical Hecke algebra}

We study the Hecke algebra $\hecke_\epsilon(\tilde T; T^\circ)$ in order to classify the spherical $\epsilon$-genuine representations of $\tilde T$.  Define $\tilde Z = Z(\tilde T)$ and define $Z^\dag = \tilde Z / \mu_n$ as before.  Define $\tilde M = Z_{\tilde T}(T^\circ)$ and let $M^\dag = \tilde M / \mu_n$.  These groups fit into a commutative diagram with exact rows:
$$\begin{tikzcd}
\mu_n \inarrow{r} \arrow{d}{=} & \tilde Z \onarrow{r} \inarrow{d} & Z^\dag \inarrow{d} \\
\mu_n \inarrow{r} & \tilde M \onarrow{r} & M^\dag.
\end{tikzcd}$$  

The main technical tool for understanding the Hecke algebra $\hecke_\epsilon(\tilde T; T^\circ)$ is the following theorem, which mostly follows from \cite[\S 6]{MWTori}.
\begin{thm}
\label{CommUnrTorus}
$\tilde M$ is a maximal commutative subgroup of $\tilde T$.  Moreover,
\begin{enumerate}
\item
$m \in M^\dag$ if and only if $\val(m) \in Y^{\sharp \Fr}$.
\item
$\tilde M = \tilde Z \cdot T^\circ$.
\end{enumerate}
\end{thm}
\proof
The fact that $m \in M^\dag$ if and only if $\val_T(m) \in Y^{\sharp \Fr}$ is the content of \cite[Proposition 6.4]{MWTori}.

Since $\tilde M$ contains $\mu_n T^\circ$, and $\tilde Z$ contains $\mu_n$, the second assertion is equivalent to the assertion that
$$\val_T(M^\dag) = \val_T(Z^\dag).$$
For this, suppose that $\val(m) \in Y^{\sharp \Fr}$.  Then \cite[Theorem 5.8]{MWTori} implies that there exists $t^\circ \in T^\circ$ such that $m t^\circ \in Z^\dag$.  In the notation of {\em loc.~cit.}, if $m = y_1(\varpi) \bar y_2(\zeta_E) t^1$, then set $t^\circ = \bar y_2(\zeta_E) t^1 \in T^\circ$ and note that $m \cdot (t^\circ)^{-1} \in Z^\dag$.  Since $\val(m) = \val(m t^\circ)$, we find that $\val_T(M^\dag) = \val_T(Z^\dag)$.  Therefore
$$\tilde M = \tilde Z \cdot T^\circ.$$

Finally, we recall from \cite[Corollary 6.5]{MWTori} that $\tilde M$ is abelian.  It is a maximal commutative subgroup of $\tilde T$, as any larger subgroup would fail to centralize $T^\circ$.
\qed

\begin{corollary}
Every $\epsilon$-genuine irreducible representation of $\tilde T$ has dimension
$$\zind(\tilde T) = \# \left( Y^{\Fr} / Y^{\sharp \Fr} \right).$$
\end{corollary}
\proof
We have found that $\tilde M$ is a maximal commutative subgroup of $\tilde T$, so every genuine irrep of $\tilde T$ is induced from a character of $\tilde M$.  The valuation $\val$ provides an isomorphism:
$$\val \From \tilde T / \tilde M \xrightarrow{\sim} Y^{\Fr} / Y^{\sharp \Fr}.$$
\qed

\begin{remark}
The formula of \cite[Proposition 4.2]{L-S} is a special case of the above.
\end{remark}

\begin{definition}
The \defined{support} of $\hecke_\epsilon(\tilde T; T^\circ)$ is
$$\{ \tilde t \in \tilde T : f(\tilde t) \neq 0 \text{ for some } f \in \hecke_\epsilon(\tilde T; T^\circ) \}.$$
\end{definition} 

Note that the Hecke algebra is supported on a set of the form $\Lambda \mu_n T^\circ$, for some subset $\Lambda \subset \tilde T / \mu_n T^\circ = T / T^\circ$.  The valuation isomorphism $\val_T \From T / T^\circ \xrightarrow{\sim} Y^{\Fr}$ suggests that we describe $\Lambda$ as a subset of $Y^{\Fr}$.

\begin{corollary}
\label{ToralSupport}
The support of $\hecke_\epsilon(\tilde T; T^\circ)$ equals $\tilde M$.  In other words, the valuation of the support of the Hecke algebra is $Y^{\sharp \Fr}$.
\end{corollary}
\proof
If $\tilde t \in \tilde T$ and $\tilde t \not \in \tilde M$, then there exists $k \in T^\circ$ such that $k \tilde t k^{-1} \tilde t^{-1} = \zeta \neq 1$.  It follows that for all $f \in \hecke_\epsilon(\tilde T; T^\circ)$,
$$f(\tilde t) = f(k \tilde t k^{-1}) = f(\zeta \tilde t) = \epsilon(\zeta) f(\tilde t)$$
and so $f(\tilde t) = 0$ (since $\epsilon$ is an injective character).   

On the other hand, if $\tilde m \in \tilde M$ then we may define a function $\delta_{\tilde m}$ supported on $\tilde m T^\circ \mu_n$, by the formula
\begin{equation}
\label{EpsilonRelation}
\delta_{\tilde m}(\tilde m k \zeta) = \epsilon(\zeta), \text{ for all } k \in T^\circ, \zeta \in \mu_n.
\end{equation}
\qed

In what follows, define $\Lambda = Y^{\sharp \Fr}$, and define $\tilde \Lambda = \tilde M / T^\circ$.  Since the valuation identifies $M^\dag / T^\circ$ with $\Lambda$, we get a commutative diagram of abelian groups with exact rows and columns as in \cite[\S 6]{MWTori}.
$$\begin{tikzcd}
\phantom{} & T^\circ \arrow{r}{=} \inarrow{d} & T^\circ \inarrow{d} \\
\mu_n \inarrow{r} \arrow{d}{=} & \tilde M \onarrow{r} \onarrow{d} & M^\dag \onarrow{d}{\val} \\
\mu_n \inarrow{r} & \tilde \Lambda \onarrow{r} & \Lambda.
\end{tikzcd}$$

The previous Corollary gives a set of vectors spanning $\hecke_\epsilon(\tilde T; T^\circ)$ as a complex vector space; for any $\tilde \lambda = \tilde m T^\circ \in \tilde \Lambda$, define $\delta_{\tilde \lambda} = \delta_{\tilde m}$.  These do not form a basis, due to the linear relations:
$$\delta_{\zeta \tilde \lambda} = \epsilon(\zeta)^{-1} \cdot \delta_{\tilde \lambda}, \text{ for all } \zeta \in \mu_n, \tilde \lambda \in \tilde \Lambda.$$

\begin{thm}
\label{ToralHecke}
The $\epsilon$-genuine spherical Hecke algebra $\hecke_\epsilon(\tilde T; T^\circ)$ is commutative.  The map $\tilde \lambda \mapsto \delta_{\tilde \lambda}$ extends to a ring isomorphism
$$\CC_\epsilon[\tilde \Lambda] \defeq \frac{ \CC[\tilde \Lambda] }{ \langle \zeta - \epsilon(\zeta)^{-1} : \zeta \in \mu_n \rangle } \xrightarrow{\sim} \hecke_\epsilon(\tilde T; T^\circ).$$
\end{thm}
\proof
Consider two elements $\tilde \lambda, \tilde \mu \in \tilde \Lambda$, representing cosets $\tilde u T^\circ$ and $\tilde v T^\circ$ (with $\tilde u, \tilde v \in \tilde M$) respectively.  The convolution of $\delta_{\tilde \lambda}$ and $\delta_{\tilde \mu}$ can be computed directly:
\begin{align*}
[\delta_{\tilde \lambda} \ast \delta_{\tilde \mu}](\tilde t) &= \int_{T} \delta_{\tilde \lambda}(h) \delta_{\tilde \mu}(h^{-1} \tilde t) dh, \\
&=  n^{-1} \sum_{\zeta \in \mu_n} \int_{\zeta \tilde u T^\circ} \epsilon(\zeta) \delta_{\tilde \mu}(\tilde h^{-1} \tilde t) d \tilde h, \\
&=   n^{-1} \sum_{\zeta \in \mu_n} \int_{T^\circ} \epsilon(\zeta) \delta_{\tilde \mu}(\zeta^{-1} \tilde u^{-1} k^{-1} \tilde t) dk, \\
&=  \int_{T^\circ} \delta_{\tilde \mu}(k^{-1} \tilde u^{-1}  \tilde t) dk, \\
&= \begin{cases} \epsilon(\zeta) & \text{ if } \tilde t \in \zeta \tilde u \tilde v T^\circ  \text{ for some } \zeta \in \mu_n; \\ 0 & \text{ otherwise.} \end{cases}
\end{align*}
Thus we find that
$$\delta_{\tilde \lambda} \ast \delta_{\tilde \mu} = \delta_{\tilde \lambda \cdot \tilde \mu}.$$
Equation \eqref{EpsilonRelation} demonstrates that the map $\tilde \lambda \mapsto \delta_{\tilde \lambda}$ is a well-defined ring homomorphism.  The fact that it is an isomorphism follows from the bijectivity of $\val_T \From T / T^\circ \To \Lambda$.  Since $\tilde \Lambda$ is commutative, we find that the Hecke algebra is also commutative.
\qed

In other words, the $\epsilon$-genuine spherical Hecke algebra is the coordinate ring of a ``commutative quantum torus'' at an $n^{\th}$ root of unity, as in \cite[Remark 6.9]{MWTori}. 

\begin{corollary}
\label{uniquespherical}
Let $(\pi, V)$ be an $\epsilon$-genuine irreducible spherical representation of $\tilde T$.  Then the space $V^{T^\circ}$ of $T^\circ$-spherical vectors is one-dimensional.
\end{corollary}

\subsection{Functorial description of the Hecke algebra}

Beginning with a central extension $\alg{K}_2 \Into \alg{T}' \Onto \alg{T}$ over $\OO$ and a positive $n$ dividing $q-1$, we have constructed an extension of abelian groups
$$\mu_n \Into \tilde \Lambda \Onto \Lambda,$$
through a multi-step process:  take $F$-points, push out via the Hilbert symbol, pull back to $\{t \in T : \val(t) \in \Lambda \}$, quotient by $T^\circ$.  The result is a commutative extension of $\Lambda$ by $\mu_n$, an object of $\Cat{Ext}(\Lambda, \mu_n)$.  

Up to natural isomorphism, the result can also be obtained by the following process:  let $\alg{T}_\Lambda$ denote the (split) torus over $\OO$ with cocharacter lattice $\Lambda = Y^{\sharp \Fr}$; pull back the central extension to obtain $\alg{K}_2 \Into \alg{T}_{\Lambda}' \Onto \alg{T}_{\Lambda}$; take $F$-points, push out via the Hilbert symbol, and quotient by $T_\Lambda^\circ = \alg{T}_{\Lambda}(\OO)$.

Thus to understand the extension $\mu_n \Into \tilde \Lambda \Onto \Lambda$, it suffices (replacing $\alg{T}$ by $\alg{T}_{\Lambda}$ if necessary) to assume that $Y = \Lambda$, i.e., that $\alg{T}$ is a split torus and $Y = Y^\sharp$.

So now assume $\alg{T}$ be a split torus over $\OO$ with character lattice $Y$.  If $Q \From Y \To \ZZ$ is a quadratic form (and the integer $n$ is fixed as always), we say that $Q$ is \defined{sharp} if
$$Y = Y^\sharp = \{ y \in Y : Q(y + y') - Q(y) - Q(y') \in n \ZZ \text{ for all } y' \in Y \}.$$

Define the Picard category $\Cat{C}(Y) = \Cat{CExt}_\OO(\alg{T}, \alg{K}_2)$.  Our conventions for Picard categories follow \cite{MWIntegral} (following \cite{SGA4T3}), and in particular all Picard categories will be strictly commutative.  Each object of this category yields (after \cite{B-D} and \cite{MWIntegral}) a ``first invariant,'' a quadratic form $Q \From Y \To \ZZ$.  Moreover, this defines a bijection from the set of isomorphism classes in $\Cat{C}(Y)$ and the set of quadratic forms $\Sym^2 X$.  Define $\Cat{C}^\sharp(Y)$ to be the full subcategory of objects whose first invariants are sharp quadratic forms.    
 
Define $\Cat{E}(Y) = \Cat{Ext}(Y, \mu_n)$, the Picard category of abelian extensions $\mu_n \Into \tilde Y \Onto Y$, with monoidal structure given by the Baer sum.  For any object $\alg{K}_2 \Into \alg{T}' \Onto \alg{T}$ of $\Cat{C}^\sharp(Y)$, define $\Fun{H}(\alg{T}') = \tilde Y$, where $\mu_n \Into \tilde Y \Onto Y$ is the extension described by the construction above (we use $\tilde Y$ instead of $\tilde \Lambda$ since $Y = \Lambda$ here).  Each step in the construction is given by an additive functor of Picard categories, and thus the construction of $\tilde Y$ defines an additive functor of Picard categories,
$$\Fun{H} \From \Cat{C}^\sharp(Y) \To \Cat{E}(Y).$$

On the other hand, the second Brylinski-Deligne invariant (adapted to the integral setting in \cite{MWIntegral}) defines a central extension of groups $\OO^\times \Into D^\circ \Onto Y$.  Consider the surjective homomorphism $h_n \From \OO^\times \Onto \mu_n$, sending $w \in \OO^\times$ (with reduction $\bar w \in \FF_q^\times$) to the Teichm\"uller lift of $\bar w^{(q-1)/n}$.  Pushing out via $h_n$ gives an extension
$$\mu_n \Into (h_n)_\ast D^\circ \Onto Y.$$
If $Q$ (the first invariant of $\alg{T}'$) is sharp, then $(h_n)_\ast D^\circ$ is commutative by the commutator formula of \cite[Proposition 3.11]{B-D}.  The construction of $(h_n)_\ast D^\circ$ from $\alg{T}'$, via the second invariant and pushing out, defines an additive functor of Picard categories
$$\Fun{D} \From \Cat{C}^\sharp(Y) \To \Cat{E}(Y).$$

As the functor $\Fun{H}$ determines the $\epsilon$-genuine spherical Hecke algebra, while $\Fun{D}$ is determined by the Brylinski-Deligne invariants alone, the following theorem describes the Hecke algebra in terms of the Brylinski-Deligne invariants.

\begin{figure}
$$\begin{tikzcd}
\phantom{a} & \alg{K}_2 \Into \alg{T}' \Onto \alg{T} \arrow{dl}[swap]{\text{$F$-points}} \arrow{dr}{\text{Brylinski-Deligne invariant}} \arrow[bend left=20]{dddl}{\Fun{H}}  \arrow[bend right=20]{dddr}[swap]{\Fun{D}} & \phantom{a} \\
\alg{K}_2(F) \Into \alg{T}'(F) \Onto T  \arrow{d}[swap]{\text{pushout via $\Hilb_n$}} & &  \OO^\times \Into D^\circ \Onto Y  \arrow{dd}{\text{pushout via $h_n$}}\\
\mu_n \Into \tilde T \Onto T  \arrow{d}[swap]{\text{quotient by $T^\circ$}} & &  \\
\mu_n \Into \tilde Y \Onto Y\arrow[Rightarrow]{rr}{N}[swap]{\sim}& &  \mu_n \Into (h_n)_\ast D^\circ \Onto Y
\end{tikzcd}$$
\caption{The two functors $\Fun{H}$ and $\Fun{D}$, beginning with an object $\alg{K}_2 \Into \alg{T}' \Onto \alg{T}$ of $\Cat{C}^\sharp(Y)$.}
\label{House}
\end{figure}

\begin{thm}
There is a natural isomorphism of additive functors $N \From \Fun{H} \xRightarrow{\sim} \Fun{D}$.
\end{thm}
\proof
We follow the methods of \cite[Theorem 4.4]{MWIntegral}, to reduce the proof to the case of ``incarnated'' extensions of tori.  Begin with an extension $\alg{K}_2 \Into \alg{T}' \Onto \alg{T}$ over $\OO$, incarnated by an element $C \in X \otimes X$ as in \cite[\S 3.10, 3.11]{B-D}.  Thus, if $C = \sum_{ij} c_{ij} x_i \otimes x_j$, then $\alg{T}' = \alg{T} \times \alg{K}_2$ (as a Zariski sheaf) with multiplication given by the rule
$$(s, \alpha) \cdot (t,\beta) = \left( s t, \alpha \beta \cdot \prod_{i,j} \{ x_i(s), x_j(t) \}^{c_{ij}} \right).$$

We apply the functor $\Fun{H}$ to $\alg{T}'$, tracing down the left side of Figure \ref{House}.  The pushout of $\alg{T}'(F)$ via $\Hilb_n$ gives $\tilde T = T \times \mu_n$, with
$$(s, \zeta_1) \cdot (t, \zeta_2) = \left( st, \zeta_1 \zeta_2 \cdot \prod_{i,j} \Hilb_n( x_i(s), x_j(t) )^{c_{ij}} \right).$$
We find $T^\circ$ embedded via $t \mapsto (t,1)$; the quotient $\tilde T / T^\circ$ forms the extension $\tilde Y = Y \times \mu_n$ with 
\begin{equation}
\label{HT}
(y_1, \zeta_1) \cdot (y_2, \zeta_2) = \left( y_1 + y_2, \zeta_1 \zeta_2 \cdot \prod_{i,j} (-1)^{c_{ij} \langle x_i, y_1 \rangle \cdot \langle x_j, y_2 \rangle \cdot (q-1)/n} \right).
\end{equation}
Indeed, if $\varpi \in F^\times$ and $\val(\varpi) = 1$ then $\Hilb_n(\varpi, \varpi) = (-1)^{(q-1)/n}$.

Now we apply the functor $\Fun{D}$ to $\alg{T}'$, tracing down the right side of Figure \ref{House}.  The second Brylinski-Deligne invariant of $\alg{T}'$ is constructed in \cite[\S 3.3.2]{MWIntegral} by first taking global sections over $\alg{G}_{\mult/\OO} = \Spec \OO[\form^{\pm 1}]$,
$$H_{\zar}^0(\alg{G}_{\mult/\OO}, \alg{K}_2) \Into H_{\zar}^0(\alg{G}_{\mult/\OO},\alg{T}') \Onto \alg{T}(\OO[\form^{\pm 1}]).$$
Pulling back via $Y \To \alg{T}(\OO[\form^{\pm 1}])$, $y \mapsto y(\form)$ and pushing out via the residue map in K-theory yields the Brylinski-Deligne invariant,
$$\OO^\times \Into D^\circ \Onto Y.$$
The residue map $\partial \From H_{\zar}^0(\alg{G}_{\mult/\OO}, \alg{K}_2) \To H_{\zar}^0(\OO, \alg{K}_1) = \OO^\times$ sends the image of $\{ \form, \form \} \in \alg{K}_2(\OO[\form^{\pm 1}])$ to $-1 \in \OO^\times$.  

We find that $D^\circ = Y \times \OO^\times$ as a set, with multiplication given by
$$(y_1, u_1) \cdot (y_2, u_2) = \left( y_1 + y_2, u_1 u_2 \cdot \prod_{i,j} (-1)^{ c_{ij} \langle x_i, y_1 \rangle \cdot \langle x_j, y_2 \rangle} \right).$$
Pushing out via $(h_n)_\ast$ yields the extension $\mu_n \Into (h_n)_\ast D^\circ \Onto Y$, with $(h_n)_\ast D^\circ = Y \times \mu_n$ as a set, and multiplication given by
\begin{equation}
\label{DT}
(y_1, \zeta_1) \cdot (y_2, \zeta_2) = \left( y_1 + y_2, \zeta_1 \zeta_2 \cdot \prod_{i,j} (-1)^{c_{ij} \langle x_i, y_1 \rangle \cdot \langle x_j, y_2 \rangle \cdot (q-1)/n} \right).
\end{equation}

The coincidence between \eqref{HT} and \eqref{DT} gives an isomorphism, 
$$N \From \Fun{H}(\alg{T}') \xrightarrow{\sim} \Fun{D}(\alg{T}').$$

The automorphism group of $\alg{T}'$ in the category $\Cat{CExt}_\OO(\alg{T}, \alg{K}_2)$ is identified with $\Hom(Y, \OO^\times) = X \otimes \OO^\times$, by the main results of \cite{MWIntegral}.  Suppose that $x \otimes w \in X \otimes \OO^\times$, and consider the corresponding automorphism $\alpha_{x \otimes w}$ of $\alg{T}'$.

On $\alg{T}' = \alg{T} \times \alg{K}_2$, this automorphism has the following effect.
$$\alpha_{x \otimes w}(t, \alpha) = (t, \alpha \cdot \{x(t), w \} ).$$
Tracing this down the left side of Figure \ref{House}, we find an automorphism of $\tilde T = T \times \mu_n$,
$$\alpha_{x \otimes w}(t, \zeta) = (t, \zeta \cdot \Hilb_n(x(t), w) ).$$
On the quotient $\tilde T / T^\circ = \tilde Y$, we find
\begin{equation}
\label{AH}
\alpha_{x \otimes w}(y, \zeta) = \left( y, \zeta \cdot \Theta(\bar w)^{ \langle x,y \rangle \cdot (q-1)/n} \right),
\end{equation}
where $\Theta(\bar w)$ denotes the Teichm\"uller lift of $\bar w \in \FF_q^\times$.  Indeed, if $\varpi \in F^\times$ and $\val(\varpi) = 1$, then $\Hilb_n(\varpi, w) = \Theta(\bar w)^{(q-1)/n}$.

Now we trace $\alpha_{x \otimes w}$ down the right side of Figure \ref{House}.  On the second Brylinski-Deligne invariant, $D^\circ = Y \times \OO^\times$, the automorphism has the effect
$$\alpha_{x \otimes u}(y,u) = (y, u \cdot  w^{\langle x,y \rangle}),$$
following \cite[\S 3.11]{B-D}.  Pushing forward via $h_n$, recalling $h_n(w) = \Theta(\bar w)^{(q-1)/n}$, we find that $\alpha_{x \otimes u}$ has the following effect on $(h_n)_\ast D^\circ = Y \times \mu_n$.
\begin{equation}
\label{AD}
\alpha_{x \otimes u}(y, \zeta) = \left( y, \zeta \cdot \Theta(\bar w)^{  \langle x,y \rangle \cdot (q-1)/n} \right).
\end{equation}
The coincidence between \eqref{AH} and \eqref{AD} implies that that isomorphism $N \From \Fun{H}(\alg{T}') \xrightarrow{\sim} \Fun{D}(\alg{T}')$ respects automorphisms.

Next, consider two elements $C_0, C \in X \otimes X$ such that $C_0(y,y) = C(y,y) = Q(y)$ for some quadratic form $Q$.  Write $\alg{T}_{C_0}'$ and $\alg{T}_C'$ for the resulting central extensions of $\alg{T}$ by $\alg{K}_2$ over $\OO$.  Let $A = C - C_0 \in X \otimes X$.  Writing $A = \sum_{i,j} a_{ij} x_i \otimes x_j$, an isomorphism from $\alg{T}_{C_0}'$ to $\alg{T}_C'$ is given by
$$\alpha(s, \kappa) = \left( s, \kappa \cdot \prod_{i,j} \{ x_i(s), x_j(s) \}^{a_{ij}} \right).$$
We apply the functor $\Fun{H}$ to the isomorphism $\alpha$, to find an isomorphism from $\tilde Y_{C_0} = Y \times \mu_n$ to $\tilde Y_C = Y \times \mu_n$,
$$\alpha(y, \zeta) = \left( y, \zeta \cdot \prod_{i,j} (-1)^{a_{ij} \langle x_i, y_1 \rangle \langle x_j, y_2 \rangle \cdot (q-1)/n } \right).$$
Applying the functor $\Fun{D}$ to the isomorphism $\alpha$ yields the same (as the extensions are identified via $N$) isomorphism from $(h_n)_\ast D_{C_0}^\circ$ to $(h_n)_\ast D_C^\circ$.  Therefore, $N$ respects the isomorphism $\alpha$.  

The groupoid $\Cat{C}^\sharp(Y)$ is equivalent to its full subcategory with objects the incarnated extensions $\{ \alg{T}_C' : C \in X \otimes X \}$.  All morphisms in this full subcategory are obtained by composing automorphisms and morphisms of the form $\alpha$ as above.  Hence there is a unique extension of $N$ to a natural isomorphism of functors,
$$N \From \Fun{H} \xRightarrow{\sim} \Fun{D}.$$

One may check that $N$ is compatible with the additive structure on these functors, since the Baer sum of incarnated extensions in $\Cat{C}^\sharp$ is obtained by adding elements of $X \otimes X$.  Hence $N$ is a natural isomorphism of additive functors, as in \cite[\S 1.4.6]{SGA4T3}
\qed

\section{Automorphic representations}
\label{Automorphic}
Now let $\alg{T}$ be a torus over a global field $F$.  Let $\alg{K}_2 \Into \alg{T}' \Onto \alg{T}$ be a central extension defined over $F$.  Let $\VV$ be the set of places of $F$, $\VV_{\fin}$ the nonarchimedean places, and $\VV_\infty$ the archimedean places.  As noted in \cite[\S 10.5]{B-D}, there exists a finite set of places $S \subset \VV$, and a model of the central extension $\alg{T}'$ defined over the $S$-integers $\OO_S \subset F$.  Any two choices of model become isomorphic after restricting to the $S'$-integers for some $S' \supset S$, and any two such isomorphisms become equality after restriction to some further $S'' \supset S'$.  In this way, we find for almost all $v \in \VV_{\fin}$, an integral model of $\alg{K}_2 \Into \alg{T}' \Onto \alg{T}$ over the local ring $\OO_v$.  See \cite{MWIntegral} for more on integral models of these central extensions.

Let $n$ be a positive integer for which $\# \mu_n(F) = n$.  We obtain a central extension
$$\mu_n \Into \tilde T_\AA \Onto T_\AA,$$
endowed with a canonical splitting over $T_F$.  The integral models provide splittings over $T_v^\circ$ for almost all $v \in \VV_{\fin}$.  

\subsection{Irreducible representations}

Our approach here to the representation theory of $\tilde T_\AA$ follows Flath \cite{Flath}.  First, let $\AA_{\fin}$ be the finite adeles.  Pulling back we find an extension
$$\mu_n \Into \tilde T_{\AA_{\fin}} \Onto T_{\AA_{\fin}},$$
and $T_{\AA_{\fin}}$ is a totally disconnected lc-group.  By construction, $\tilde T_{\AA_{\fin}}$ arises as a quotient of a restricted direct product in the following way:  define
$$\dtilde T_{\AA_{\fin}} = \biprod_{v \in \VV_{\fin}} \left( \tilde T_v, T_v^\circ \right)$$
Here the notation refers to the restricted direct product of the covers $\tilde T_v$, with respect to the compact open subgroups $T_v^\circ$ (defined almost everywhere).  This restricted direct product fits into a short exact sequence, and $\tilde T_{\AA_{\fin}}$ is the pushout via the product map below.
$$\begin{tikzcd}
\bigoplus_{v \in \VV_{\fin}} \mu_n \inarrow{r} \onarrow{d}{\Pi} & \dtilde T_{\AA_{\fin}} \onarrow{r} \onarrow{d} & T_{\AA_{\fin}} \arrow{d}{=}  \\
\mu_n \inarrow{r} & \tilde T_{\AA_{\fin}} \onarrow{r} & T_{\AA_{\fin}}
\end{tikzcd}$$

Following \cite[Example 2]{Flath}, the restricted direct product $\dtilde{T}_{\AA_{\fin}}$ is a totally disconnected lc-group, and its Hecke algebra is isomorphic (as an idempotented $\CC$-algebra) to the tensor product of the Hecke algebras $\hecke(\tilde T_v)$ of locally constant compactly supported measures, with respect to their idempotents $\Char(T_v^\circ)$ (the characteristic function of $T_v^\circ$) at almost all places.  The commutativity of the Hecke algebras from Theorem \ref{ToralHecke} implies a factorization of irreducible admissible representations.
\begin{proposition}
Every irreducible admissible representation $\pi_{\fin}$ of $\dtilde T_{\AA_{\fin}}$ is factorizable:  there exists a unique family of equivalence classes $([\pi_v])_{v \in \VV_{\fin}}$ of irreducible admissible representations of each $\tilde T_v$, such that $\pi_v$ is $T_v^\circ$-spherical for almost all $v \in \VV_{\fin}$ and $\pi_{\fin}$ is isomorphic to the restricted tensor product of representations $\bigotimes_{v \in \VV_{\fin}} \pi_v$ with respect to some choice of nonzero spherical vectors at almost all places.
\end{proposition}
\proof
This is a direct result of \cite[Theorem 2]{Flath}
\qed

To include the archimedean places, define
$$\dtilde T_\AA = \dtilde T_{\AA_{\fin}} \times \prod_{v \in \VV_\infty} \tilde T_v.$$
Then the extension $\tilde T_\AA$ is the pushout via the product map below.
$$\begin{tikzcd}
\bigoplus_{v \in \VV} \mu_n \inarrow{r} \onarrow{d}{\Pi} & \dtilde T_{\AA} \onarrow{r} \onarrow{d} & T_{\AA} \arrow{d}{=}  \\
\mu_n \inarrow{r} & \tilde T_\AA \onarrow{r} & T_\AA
\end{tikzcd}$$

When $v \in \VV_\infty$, the admissible representations of $\tilde T_v$ (i.e., $(\Lie{t}, \tilde T^\circ)$-modules) are admissible modules over the Hecke algebra $\hecke(\Lie{t}, \tilde T^\circ)$ of left and right $\tilde T_v^\circ$-finite distributions on $\tilde T_v$ with support in $\tilde T_v^\circ$.  In this way, Flath includes archimedean places in the factorization of admissible representations as well.  Here, an \defined{admissible representation} of $\dtilde T_\AA$ will mean a representation which is simultaneously an admissible $\dtilde T_{\AA_{\fin}}$ representation and an admissible representation of $\tilde T_v$ for all $v \in \VV_\infty$.

\begin{proposition}
Every irreducible admissible representation $\pi$ of $\dtilde T_{\AA}$ is factorizable:  there exists a unique family of equivalence classes $([\pi_v])_{v \in \VV}$ of irreducible admissible representations of each $\tilde T_v$, such that $\pi_v$ is $T_v^\circ$-spherical for almost all $v \in \VV_{\fin}$ and $\pi$ is isomorphic to the restricted tensor product of representations $\bigotimes_{v \in \VV} \pi_v$ with respect to some choice of nonzero spherical vectors at almost all places $v \in \VV_{\fin}$.
\end{proposition}

This result descends to a result for the cover $\tilde T_\AA$.  An \defined{admissible representation} of $\tilde T_\AA$ will mean a representation which becomes admissible after pulling back to $\dtilde T_\AA$.
\begin{corollary}
Every irreducible admissible $\epsilon$-genuine representation $\pi$ of $\tilde T_{\AA}$ is factorizable:  there exists a unique family of equivalence classes $([\pi_v])_{v \in \VV}$ of irreducible $\epsilon$-genuine admissible representations of each $\tilde T_v$, such that $\pi_v$ is $T_v^\circ$-spherical for almost all $v \in \VV_{\fin}$ and $\pi$ is isomorphic to the restricted tensor product of representations $\bigotimes_{v \in \VV} \pi_v$ with respect to some choice of nonzero spherical vectors at almost all places $v \in \VV_{\fin}$.
\end{corollary}
\proof
The representations of $\tilde T_\AA$ are those of $\dtilde T_\AA$ on which
$$\Ker \left( \bigoplus_{v \in \VV} \mu_n \Onto \mu_n \right)$$
acts trivially.  Therefore the $\epsilon$-genuine representations of $\tilde T_\AA$ are those representations of $\dtilde T_\AA$ on which each $\mu_n$ summand acts via $\epsilon$.
\qed

The unitary version of this factorization is similar.  This is essentially the content of \cite[Theorem 4]{Flath}, which relies on \cite{Bern} for admissibility of irreducible unitary representations of $p$-adic groups and \cite{GGPS}.  The factorization also follows from early work of Moore \cite{Moo} as follows.  

A {\em unitary} representation of $\tilde T_\AA$ (or $\dtilde T_\AA$) will always mean a continuous unitary representation on a Hilbert space.  If $\pi$ is an irreducible unitary representation of $\dtilde T_\AA$, then for any place $v \in \VV$ we get a factorization $\pi \isom \pi_v \hat \otimes \pi_{v'}$, where $\pi_{v'}$ is an irreducible unitary representation of $\dtilde T_{\AA_{v'}} = \biprod_{u \in \VV - v} (\tilde T_u, T_u^\circ)$.  The equivalence class of each $\pi_v$ is uniquely determined by that of $\pi$.  In \cite[Lemma 6.3]{Moo}, Moore demonstrates that $\pi_v$ is $T_v^\circ$-spherical for almost all $v \in \VV$.  Commutativity of the spherical Hecke algebras and \cite[Theorems 6,7]{Moo}, imply that the factorization above gives a bijection between the sets of...
\begin{itemize}
\item
equivalence classes of unitary irreducible representations $[\pi]$ of $\dtilde{T}_\AA$;
\item
families $([\pi_v])_{v \in \VV}$ of equivalence classes of unitary irreducible representations of $\tilde T_v$, which are $T_v^\circ$-spherical for almost all $v \in \VV$.
\end{itemize}

In our setting, where all irreducible admissible representations of $\tilde T_v$ are finite-dimensional, we get a simple version of \cite[Theorem 4]{Flath}.
\begin{proposition}
Suppose that $\pi$ is an irreducible unitary $\epsilon$-genuine representation of $\tilde T_\AA$.  Then there exist, unique up to equivalence, irreducible unitary $\epsilon$-genuine representatiosn $\pi_v$ of $T_v$ for all $v \in \VV$ such that $\pi \isom \widehat{\bigotimes} \pi_v$ (the Hilbert space restricted tensor product with respect to normalized spherical vectors).
\end{proposition}
Here, the irreducible unitary $\epsilon$-genuine representations of each $\tilde T_v$ are all finite-dimensional admissible representations; there is no need to take $T_v^\circ$-finite vectors.

The factorization of $\epsilon$-genuine irreducible representations of $\tilde T_\AA$, together with the classification of $\epsilon$-genuine irreducible representations of $\tilde T_v$ by their central character, gives the following result.
\begin{thm}
\label{SvNGlobal}
Factorization and taking central characters gives a bijection between the sets of...
\begin{itemize}
\item
equivalence classes $[\pi]$ of irreducible admissible $\epsilon$-genuine representations of $\tilde T_\AA$;
\item
families $([\pi_v])_{v \in \VV}$ of equivalence classes of irreducible admissible $\epsilon$-genuine representations of each $\tilde T_v$, almost all of which are spherical;
\item
families $(\chi_v)_{v \in \VV}$ of elements $\chi_v \in \Hom_\epsilon(Z(\tilde T_v), \CC^\times)$, almost all of which are trivial on $Z(\tilde T_v) \cap T_v^\circ$;
\item
$\epsilon$-genuine characters $\chi \in \Hom_\epsilon(Z(\tilde T_\AA), \CC^\times)$.
\end{itemize}
These bijections preserve the subsets of unitary representations throughout.
\end{thm}
\proof
The bijection between the first two sets is the factorization result above.  The bijection between the second and third set is the classification of irreducible admissible representations of $\tilde T_v$ for each local field, Theorem \ref{LocalSvN}.  The only thing to observe is that spherical irreps of $\tilde T_v$ are those whose central character is trivial on $Z(\tilde T_v) \cap T_v^\circ$; indeed such a central character extends to a character of the maximal abelian subgroup $Z(\tilde T_v) \cdot T_v^\circ$ which is trivial on $T_v^\circ$, whose induction is a spherical irreducible representation of $\tilde T_v$.  

The bijection between the third set and the fourth follows from the description of the center of $\tilde T_\AA$; recall that $Z_\AA^\dag$ is the restricted direct product of $Z_v^\dag$ with respect to their open compact subgroups $Z_v^\circ = Z_v^\dag \cap T_v^\circ$.  The bijection follows.
\qed

The previous theorem classifies the $\epsilon$-genuine unitary irreducible representations of $\tilde T_\AA$, by using the local classification (Theorem \ref{LocalSvN}) and factorization.  But we also require a purely global understanding of the unitary irreducible representations of $\tilde T_\AA$; this requires a version of the Stone von-Neumann theorem applicable in a setting of locally compact, not necessarily separable, topological groups.  A similar result (for similar goals) is given in \cite[\S 0.3]{K-P}, but their treatment leaves much to the reader.  A modern treatment, in sufficient generality, can be found in work of L\"udeking and Poguntke \cite[Theorem 2.3]{L-P}.  A direct consequence of their work is the following theorem.
\begin{thm}
\label{GlobalPolarizedSvN}
Let $\tilde M$ be a maximal abelian subgroup of $\tilde T_\AA$, and let $\chi$ be an $\epsilon$-genuine unitary character of $\tilde Z_\AA$.  Suppose that $\tilde Z_\AA \backslash \tilde M$ is discrete and $\tilde T_\AA / \tilde M$ is compact.  Let $\pi_\chi$ be the (unique up to equivalence) unitary irreducible representations of $\tilde T_\AA$ with central character $\chi$.  There exists an extension of $\chi$ to a unitary character $\chi_M$ of $\tilde M$.  For any such extension $\chi_M$, the unitary induced representation (in the sense of Blattner \cite{Bla}) $\Ind_{\tilde M}^{\tilde T_\AA} \chi_M$ is irreducible and equivalent to $\pi_\chi$.
\end{thm}
\proof
Let $Q = \tilde Z_\AA \backslash \tilde T_\AA$ (a locally compact abelian group, with the quotient topology), and consider the nondegenerate continuous alternating form $\gamma \From Q \times Q \To \mu_n$,
$$\gamma(\bar t_1, \bar t_2) = \epsilon \left( \Comm(t_1, t_2) \right).$$

Let $M = \tilde Z_\AA \backslash \tilde M$ be the image of the maximal abelian subgroup $\tilde M$ in $Q$.  This is a closed subgroup of $Q$, isotropic for $\gamma$.  Since $\tilde M$ is maximal among abelian subgroups of $\tilde T_\AA$, we find that
$$M^\perp = \{ q \in Q : \gamma(m,q) = 1 \text{ for all } m \in M \} = M.$$
Nondegeneracy of $\gamma$ implies that the continuous map $M \To (Q/M)^\ast$, $m \mapsto \gamma(m, \bullet)$ is injective with dense image.

Since $M$ is discrete, and $(Q/M) \isom \tilde T_\AA / \tilde M$ is compact, we find that $(Q/M)^\ast$ is discrete; it follows that the injective map $M \To (Q/M)^\ast$ with dense image must be a topological isomorphism.  In other words, $M$ is a {\em quasi-polarization} of the {\em quasi-symplectic space} $Q$, in the terminology of \cite[Def 1.1]{L-P}.  

The theorem now follows from Theorem \ref{SvNGlobal}, and from \cite[Theorem 2.3, Corollary 2.6]{L-P}.
\qed

\subsection{Definition of automorphic representations}

\begin{definition}
An \defined{$\epsilon$-genuine automorphic central character} ($\epsilon$-genuine \defined{ACC}, for short) for $\tilde T_\AA$ is a continuous homomorphism $\chi \From \tilde Z_\AA \To \CC^\times$ such that $\chi$ is trivial on $T_F \cap \tilde Z_\AA$ and the restriction of $\chi$ to $\mu_n$ coincides with $\epsilon$.
\end{definition}

\begin{proposition}
\label{ExtendACC}
Every $\epsilon$-genuine ACC $\chi \From \tilde Z_\AA \To \CC^\times$ extends uniquely to a continuous character $T_F \tilde Z_\AA \To \CC^\times$ which is trivial on $T_F$.
\end{proposition}
\proof
Since $T_F$ is a discrete subgroup of $\tilde T_\AA$, the group $\tilde Z_\AA$ is open in $T_F \tilde Z_\AA$.  Thus the inclusion $Z_\AA \Into T_F Z_\AA$ induces a topological isomorphism of lca-groups,
$$\frac{\tilde Z_\AA}{T_F \cap \tilde Z_\AA} \xrightarrow{\sim} \frac{T_F \cdot \tilde Z_\AA}{T_F}.$$
The result follows immediately.
\qed

Let $\hecke_{\fin} = \hecke(\dtilde T_{\AA_{\fin}})$ denote the Hecke algebra of locally constant, compacly supported measures on the totally disconnected lc-group $\dtilde T_{\AA_{\fin}}$.  Let $\hecke_\infty = \bigotimes_{v \in \VV_\infty} \hecke(\Lie{t}_v, \tilde T_v^\circ)$.  Functions on $\tilde T_\AA$ will be viewed as functions on $\dtilde T_\AA$ by pulling back, as needed.  The following definition is lifted from Borel and Jacquet \cite[\S 4.2, 5.8]{BorelJacquet}, with adaptations based on the equivalences of \cite[Proposition 4.5, Corollary 5.7]{BorelJacquet}.  In particular, when $n = 1$ so that $\tilde T_\AA = T_\AA$, the definition coincides with the ``usual'' definition of automorphic forms given in \loccit.
\begin{definition}
A function $f \From \tilde T_\AA \To \CC$ is an $\epsilon$-genuine \defined{automorphic form} if it satisfies the following conditions.
\begin{enumerate}
\item
$f(\zeta \gamma \tilde t) = \epsilon(\zeta) \cdot f(\tilde t)$ for all $\gamma \in T_F$, $\zeta \in \mu_n$, and $\tilde t \in \tilde T_\AA$.
\item
$f$ is locally constant on $\tilde T_{\AA_{\fin}}$ and smooth on $\tilde T_v$ for all $v \in \VV_\infty$.
\item
There is an element $\xi = \xi_{\fin} \otimes \xi_\infty \in \hecke_{\fin} \otimes \hecke_\infty$ such that $f \ast \xi = f$.
\item
For every $y \in \tilde T_{\AA}$ and every $v \in \VV_\infty$, the function $\tilde T_v \To \CC$, $x \mapsto f(x \cdot y)$ is slowly increasing.
\item
The span of $\{ f \ast \xi : \xi \in \hecke_{\fin} \otimes \hecke_\infty \}$ is an admissible representation of $\tilde T_\AA$.
\end{enumerate}
The space of such automorphic forms will be written $\AF_\epsilon(\tilde T_\AA)$.  The subspace corresponding to a specific $\xi$ will be denoted $\AF_\epsilon(\tilde T_\AA, \xi)$.
\end{definition}

\begin{remark}
Proposition 4.5 and Corollary 5.7 of \cite{BorelJacquet} probably generalize to the cover $\tilde T_\AA$, so that to check admissibility in condition (5), it should suffice to check admissibility at a nonempty set of places including all archimedean places.
\end{remark}

\begin{definition}
An $\epsilon$-genuine \defined{automorphic representation} of $\tilde T_\AA$ is an irreducible admissible subquotient of $\AF_\epsilon(\tilde T_\AA)$.
\end{definition}

Every $\epsilon$-genuine irreducible admissible representation of $\tilde T_\AA$ has a central character $\chi \From \tilde Z_\AA \To \CC^\times$; this central character must be $\epsilon$-genuine automorphic (an ACC).  Conversely, suppose that $\chi \From \tilde Z_\AA \To \CC^\times$ is an $\epsilon$-genuine ACC.  We may consider the subspace $\AF_\chi(\tilde T_\AA)$ consisting of automorphic forms which satisfy
$$f(\tilde z \tilde t) = \chi(\tilde z) \cdot f(\tilde t), \text{ for all } \tilde z \in \tilde Z_\AA, \tilde t \in t \in \tilde T_\AA.$$

The following lemma explains why the study of $\epsilon$-genuine automorphic forms with fixed ACC is easy, from the spectral point of view.
\begin{lemma}
\label{QuotCompact}
The space $T_F \tilde Z_\AA \backslash \tilde T_\AA$ is compact.
\end{lemma}
\proof
The space $T_F \tilde Z_\AA \backslash \tilde T_\AA$ can be identified with $T_F Z_\AA^\dag \backslash T_\AA$.  It is a classical fact that $T_F T_\RR \backslash T_\AA$ is compact; since $Z_\RR^\dag$ has finite index in $T_\RR$, it follows that $T_F Z_\RR^\dag \backslash T_\AA$ is compact.  Thus its quotient  $T_F Z_\AA^\dag \backslash T_\AA$ is compact.
\qed

The spectral theory of automorphic forms on covering groups is studied by Moeglin and Waldspurger \cite{MW}.  Their work implies the following.
\begin{proposition}
Let $\chi$ be an $\epsilon$-genuine ACC and let $[\pi_\chi]$ denote the unique equivalence class of admissible representations of $\tilde T_\AA$ with central character $\chi$.  Then $\AF_\chi(\tilde T_\AA)$ is equivalent to a finite direct sum of copies of $\pi_\chi$, as an admissible representation of $\tilde T_\AA$.
\end{proposition}
\proof
This falls into the ``$P = G$'' case of cuspforms, discussed in \cite[\S 1.2.18]{MW}; from their work it follows that $\AF_\chi(\tilde T_\AA)$ is semisimple as a representation of $\tilde T_\AA$, and that for each $\xi = \xi_{\fin} \otimes \xi_\infty$, the space $\AF_\chi(\tilde T_\AA, \xi)$ is finite-dimensional.  (A priori, one should cut down further by an ideal in the universal enveloping algebra of the Lie algebra of $\prod_{v \in \VV_\infty} \tilde T_v$.  But the central character takes care of this already.)  

Since irreducible admissible representations of $\tilde T_\AA$ are determined up to isomorphism by their central character (Theorem \ref{SvNGlobal}), the space $\AF_\chi(\tilde T_\AA)$ is isotypic as a representation of $\tilde T_\AA$.  Hence the finite-dimensionality of its subspaces $\AF_\chi(\tilde T_\AA, \xi)$ implies that $\AF_\chi(\tilde T_\AA)$ is equivalent to a finite direct sum of the unique irreducible representation of $\tilde T_\AA$ with central character $\chi$,
$$\left[ \AF_\chi(\tilde T_\AA) \right] = m_\pi \cdot [\pi_\chi], \text{ for some } 0 \leq m_\pi < \infty$$ 
\qed

It is also important to study automorphic forms and representations of $\tilde T_\AA$ in the unitary setting.  For this, fix a {\em unitary} $\epsilon$-genuine ACC $\chi$.  Define $L_\chi^2(T_F \backslash \tilde T_\AA)$ to be the unitary representation of $\tilde T_\AA$, induced from the unitary character $\chi$ of $T_F \cdot \tilde Z_\AA$ (via Proposition \ref{ExtendACC}).  This coincides (see \cite[\S I.2.18]{MW}) with the Hilbert space completion of the space $ \AF_\chi(\tilde T_\AA)$ with respect to the inner product
$$\langle f_1, f_2 \rangle = \int_{T_F \tilde Z_\AA \backslash \tilde T_\AA} \overline{f_1(x)} f_2(x) dx,$$
noting that the quotient is compact (Lemma \ref{QuotCompact}) and the integrand well-defined.

The previous result gives a decomposition of unitary representations of $\tilde T_\AA$,
$$\left[ L_\chi^{2}(T_F \backslash \tilde T_\AA) \right] = m_\chi \cdot [\hat \pi_\chi], \text{ for some } 0 \leq m_\chi < \infty.$$
Here $[\hat \pi_\chi]$ denotes the equivalence class of the Hilbert space completion of $\pi_\chi$.

The description of $L_\chi^2(T_F \backslash \tilde T_\AA)$ as an induced representation, and induction in stages, gives the following multiplicity formula.
\begin{proposition}
\label{MultiplicityIndex}
Let $\tilde M$ be a maximal abelian subgroup of $\tilde T_\AA$ containing $T_F$, and assume that $\tilde Z_\AA \backslash \tilde M$ is discrete and $\tilde T_\AA / \tilde M$ is compact.  Then for any $\epsilon$-genuine unitary ACC $\chi$, we have
$$1 \leq m_{[\chi]} = \# \frac{\tilde M}{T_F \tilde Z_\AA} < \infty.$$
\end{proposition}
\proof
Pontrjagin duality gives a short exact sequence
$$(T_F \tilde Z_\AA \backslash \tilde M)^\ast \Into (T_F \backslash \tilde M)^\ast \Onto (T_F \backslash T_F \tilde Z_\AA)^\ast.$$
Let $\chi_M$ be a unitary character of $T_F \backslash \tilde M$ which extends $\chi$.  Since $\tilde Z_\AA \backslash \tilde M$ is discrete, 
$$\Ind_{T_F \tilde Z_\AA}^{\tilde M} \chi \isom \bigoplus_{\psi \in (T_F \tilde Z_\AA \backslash \tilde M)^\ast} \psi \cdot \chi_M.$$
Note that compactness of $T_F \tilde Z_\AA \backslash \tilde T_\AA$ combined with discreteness of $\tilde Z_\AA \backslash \tilde M$ implies that $T_F \tilde Z_\AA \backslash \tilde M$ is finite.

Now, inducing from $\tilde M$ to $\tilde T_\AA$, applying Theorem \ref{GlobalPolarizedSvN}, we find
$$L_\chi^2(T_F \backslash \tilde T_\AA) \isom \Ind_{\tilde M}^{\tilde T_\AA} \Ind_{T_F \tilde Z_\AA}^{\tilde M} \chi \isom \bigoplus_{\psi \in (T_F \tilde Z_\AA \backslash \tilde M)^\ast} \pi_\chi.$$
The result follows immediately.
\qed

\subsection{A multiplicity formula}

Fix an $\epsilon$-genuine \textbf{unitary} ACC $\chi$ here.  To understand the multiplicity $m_\chi$ of $\pi_\chi$ in the isotypic representation $L_\chi^2(T_F \backslash \tilde T_\AA)$, we must understand the index of $T_F \tilde Z_\AA$ in a maximal abelian subgroup of $\tilde T_\AA$.  For this, we apply Poitou-Tate global duality, for the finite Galois modules in duality:
$$\alg{\hat T}_{[n]} = \sheaf{X} \otimes \mu_n, \quad \alg{T}_{[n]} = \sheaf{Y} \otimes \mu_n.$$
Our treatment of Galois cohomology follows \cite[II.\S 6]{Serre}; we use Tate's name in notation for duality pairings, though we should acknowledge Nakayama and Poitou as well.  For any nonarchimedean place $v$, write $H_{\unr}^1(F_v, \alg{\hat T}_{[n]})$ and $H_{\unr}^1(F_v, \alg{T}_{[n]})$ for the unramified part of the cohomology groups $H_{\et}^1(F_v, \alg{\hat T}_{[n]})$ and $H_{\et}^1(F_v, \alg{T}_{[n]})$.  Define the restricted products:
$$P^1(\alg{T}_{[n]}) = \biprod_{v} \left( H_{\et}^1(F_v, \alg{T}_{[n]}), H_{\unr}^1(F_v,\alg{T}_{[n]}) \right);$$
$$P^1(\alg{\hat T}_{[n]}) = \biprod_{v} \left( H_{\et}^1(F_v, \alg{\hat T}_{[n]}), H_{\unr}^1(F_v, \alg{\hat T}_{[n]}) \right).$$

In the local pairing, $\TD_v \From H_{\et}^1(F_v, \alg{\hat T}_{[n]}) \times H_{\et}^1(F_v, \alg{T}_{[n]}) \To \mu_n$ (see \eqref{LocalTatePairing}), the unramified parts of cohomology are annihilators of each other.  The global Tate pairing is thus defined via a product of local Tate pairings.
$$\TD \From P^1(\alg{\hat T}_{[n]}) \times P^1(\alg{T}{[n]})  \xrightarrow{\prod_v \TD_v} \mu_n.$$

Within the 9-term exact sequence of Poitou-Tate global duality, three terms are:
\begin{equation}
\label{PoitouTate}
H_{\et}^1(F, \alg{\hat T}_{[n]}) \xrightarrow{\loc} P^1(\alg{\hat T}_{[n]} ) \xrightarrow{\tau} H_{\et}^1(F, \alg{T}{[n]})^\ast.
\end{equation}
Hereafter, $\loc$ stands for the natural maps from Galois cohomology over $F$ to the restricted direct product of cohomologies over each $F_v$.  The map $\tau$ is given by
$$\tau(p)(\eta) =\TD(p, \loc \eta), \text{ for all } p \in P^1(\alg{\hat T}_{[n]}), \eta \in H_{\et}^1(F, \alg{T}_{[n]}).$$
For the following theorem, we use the group $\alg{R}$ from Equation \eqref{DefOfR}.
\begin{thm}
There exists a maximal abelian subgroup $\tilde M \subset \tilde T_\AA$ such that $T_F \tilde Z_\AA \subset \tilde M$ and
$$\# \frac{\tilde M}{T_F \tilde Z_\AA} \leq \# \Ker \left( \Sha^1(\alg{R}) \To \Sha^1(\alg{T} \times \alg{\hat T}) \right) \cdot \# \Sha^1(\alg{T} \times \alg{\hat T})_{[n]}.$$
In particular, if $\Sha^1(\alg{T})_{[n]} = \Sha^1(\alg{\hat T})_{[n]} = \Sha^1(\alg{R}) = 0$, (e.g., if $\alg{T}$ is a split torus) then $T_F \tilde Z_\AA$ is a maximal abelian subgroup of $\tilde T_\AA$.
\end{thm}
\proof
On the one hand, $T_F \cdot \tilde Z_\AA$ is abelian, since $T_F$ is an abelian subgroup of $\tilde T_\AA$ (via the canonical splitting of the cover over $T_F$).  To demonstrate the theorem, begin with the centralizer $\tilde N = Z_{\tilde T_\AA}(T_F)$ and its image $N^\dag \subset T_\AA$.  Note that $N^\dag$ contains $T_F Z_\AA^\dag$.  

If $t \in T_\AA$, then the local Kummer coboundaries yield an element $\Kum t \in P^1(\alg{T}_{[n]})$; similarly, if $s \in T_F$ then the global Kummer coboundary yields $\Kum s \in H_{\et}^1(F, \alg{T}_{[n]})$.  The commutator $\Comm(t,s)$ is given by $\Comm(t,s) = \TD( \delta j \Kum t, \loc  \Kum s)$.  Therefore,
\begin{align*}
N^\dag &= \{ t \in T_\AA : \Comm(t,s) = 1 \text{ for all } s \in T_F \}, \\
&= \{ t \in T_\AA : \TD(\delta j \Kum t, \loc \Kum s) = 1 \text{ for all } s \in T_F \}.
\end{align*}
Define a subgroup of $N^\dag$,
$$E^\dag = \{ t \in T_\AA : \TD(\delta j \Kum t, \loc \eta) = 1 \text{ for all } \eta \in H_{\et}^1(F, \alg{T}_{[n]}) \}.$$
The group $Z_\AA^\dag$ has such a characterization,
$$Z_\AA^\dag = \{ t \in T_\AA : \TD(\delta j \Kum t, \Kum s) = 1 \text{ for all } s \in T_\AA \}.$$

The Kummer sequences and localization maps fit into a commutative diagram.
\begin{equation}
\label{CDForT}
\begin{tikzcd}
\phantom{a} & \Sha^1(\alg{T}_{[n]}) \arrow{r}  \inarrow{d} & \Sha^1(\alg{T})_{[n]} \inarrow{d} \\
T_{F/n} \inarrow{r}{\Kum} \arrow{d}{\loc} & H_{\et}^1(F, \alg{T}_{[n]}) \onarrow{r}{\lambda} \arrow{d}{\loc} & H_{\et}^1(F, \alg{T})_{[n]} \arrow{d}{\loc} \\
T_{\AA/n} \inarrow{r}{\Kum} & P^1(\alg{T}_{[n]}) \arrow{r}{\lambda} &P^1(\alg{T})_{[n]}
\end{tikzcd}
\end{equation}

Consider the homomorphism $\phi \From T_\AA \To \Hom(P^1(\alg{T}_{[n]}), \mu_n)$ given by
$$\phi(t) = \TD( \delta j \Kum t, \bullet).$$
Then $\phi$ restricts to an injective homomorphism
$$\bar \phi \From \frac{N^\dag}{E^\dag Z_\AA^\dag} \Into \Hom \left( \frac{\loc H_{\et}^1(F, \alg{T}_{[n]}) \cap \Kum (T_\AA) } {\loc \Kum (T_F) }, \mu_n \right).$$
A short diagram chase gives an embedding
$$\frac{\loc H_{\et}^1(F, \alg{T}_{[n]}) \cap \Kum (T_\AA) } {\loc \Kum (T_F) } \hookrightarrow \frac{\Sha^1(F, \alg{T})_{[n]}}{\Im(\Sha^1(\alg{T}_{[n]}))}.$$
We find immediately that
\begin{equation}
\label{EstimateNE}
\# \frac{N^\dag}{E^\dag Z_\AA^\dag} \leq \# \frac{ \Sha^1(F, \alg{T})_{[n]}}{\Im(\Sha^1(\alg{T}_{[n]}))}.
\end{equation}

Now Poitou-Tate duality (\ref{PoitouTate}) characterizes $E^\dag$ in the following way:
$$E^\dag = \left\{ t \in T_\AA : \delta j \Kum t \in \Im \left( H_{\et}^1(F, \alg{\hat T}_{[n]}) \xrightarrow{\loc} P^1(\alg{\hat T}_{[n]}) \right) \right\}.$$

Consider the commutative diagram below, the analogue of \eqref{CDForT} for $\alg{\hat T}$.
\begin{equation}
\label{CDForTHat}
\begin{tikzcd}
\phantom{a} & \Sha^1(\alg{\hat T}_{[n]}) \inarrow{d} \arrow{r} & \Sha^1(\alg{\hat T})_{[n]} \inarrow{d} \\
\hat T_{F/n} \inarrow{r}{\Kum} \arrow{d}{\loc} & H_{\et}^1(F, \alg{\hat T}_{[n]}) \onarrow{r}{\hat \lambda} \arrow{d}{\loc} & H_{\et}^1(F, \alg{\hat T})_{[n]} \arrow{d}{\loc} \\
\hat T_{\AA / n} \inarrow{r}{\Kum} & P^1(\alg{\hat T}_{[n]}) \arrow{r}{\hat \lambda} & P^1(\alg{\hat T})_{[n]}
\end{tikzcd}
\end{equation}

If $t \in E^\dag$, then there exists $\eta \in H_{\et}^1(F, \alg{\hat T}_{[n]})$ such that $\loc \eta = \Kum \delta j t$.  The element $\hat \lambda(\eta) \in H_{\et}^1(F, \alg{\hat T})_{[n]}$ is contained in $\Sha^1(\alg{\hat T})_{[n]}$ by a quick diagram chase.  As $t$ defines $\eta$ uniquely up to $\Sha^1(\alg{\hat T}_{[n]})$, this defines a homomorphism
$$\phi \From E^\dag \To \frac{ \Sha^1(\alg{\hat T})_{[n]}}{\Im( \Sha^1(\alg{\hat T}_{[n]}))}, \quad \phi(t) = \hat \lambda(\loc^{-1}(t)).$$

If $\phi(t) = 0$, then $\hat \lambda(\eta) \in \Im(\Sha^1(\alg{\hat T}_{[n]}))$ and therefore
$$\eta \in \Kum(\hat T_{F/n}) \cdot \Sha^1(\alg{\hat T}_{[n]}).$$
Applying $\loc$, we find that
$$\Kum \delta j t = \loc \eta = \loc \Kum(\hat t), \text{ for some } \hat t \in \hat T_{F/n}.$$
Define a subgroup of $E^\dag$,
$$H^\dag = \left\{ t \in T_\AA : \delta j \Kum t \in \Im \left( \hat T_F \xrightarrow{\loc \Kum} P^1(\alg{\hat T}_{[n]}) \right) \right\}.$$
We find an inclusion
\begin{equation}
\label{EstimateEH}
\frac{E^\dag}{H^\dag} \xhookrightarrow{\phi} \frac{ \Sha^1(\alg{\hat T})_{[n]}}{\Im( \Sha^1(\alg{\hat T}_{[n]}))}.
\end{equation}

Now we recharacterize $H^\dag$ to relate it to $Z_\AA^\dag$.  Consider the map $g \From \alg{T} \times \alg{\hat T} \To \alg{\hat T}$, given on points by:
$$g(u, \hat u) = \delta j(u) \cdot \hat u^{-n}.$$
The kernel of this map is precisely the group $\alg{R}$ (recall \eqref{DefOfR}), leading to a short exact sequence of algebraic groups over $F$:
$$\alg{R} \Into \alg{T} \times \alg{\hat T} \xtwoheadrightarrow{g} \alg{\hat T}.$$
Let $\sigma$ denote the connecting map in Galois cohomology:
$$T_F \times \hat T_F \xrightarrow{g} \hat T_F \xrightarrow{\sigma} H_{\et}^1(F, \alg{R}).$$

We use the fact that $\delta j \circ \Kum = \Kum \circ \delta j$ and $\loc \circ \Kum = \Kum \circ \loc$ (functoriality of the Kummer sequence) to find:
\begin{align*}
H^\dag &=  \left\{ t \in T_\AA : \delta j \Kum t = \loc \Kum \hat t \text{ for some } \hat t \in \hat T_F \right\} \\
&=  \left\{ t \in T_\AA : \Kum \left( \delta j t - \loc \hat t \right) = 0 \text{ for some } \hat t \in \hat T_F \right\} \\
&=  \left\{ t \in T_\AA : \left( \delta j t - \loc \hat t \right) \in n \cdot \hat T_\AA \text{ for some } \hat t \in \hat T_F \right\}.
\end{align*}
Now, if $\delta j t - \loc \hat t \in n \cdot \hat T_\AA$, then $\loc(\hat t) \in g(T_\AA \times \hat T_\AA)$.  Consider the following commutative diagram with exact rows and columns.
$$\begin{tikzcd}
\phantom{a} & \phantom{a} & \Sha^1(\alg{R}) \inarrow{d} \arrow{r} & \Sha^1(\alg{T} \times \alg{\hat T}) \inarrow{d} \\
T_F \times \hat T_F \arrow{r}{g} \inarrow{d}{\loc} & \hat T_F \arrow{r}{\sigma} \inarrow{d}{\loc} & H_{\et}^1(F, \alg{R}) \arrow{d}{\loc} \arrow{r} & H_{\et}^1(F, \alg{T} \times \alg{\hat T}) \arrow{d}{\loc} \\
T_\AA \times \hat T_\AA \arrow{r}{g} & \hat T_\AA \arrow{r}{\sigma} & P^1(\alg{R}) \arrow{r} & P^1(\alg{T} \times \alg{\hat T})
\end{tikzcd}$$

The condition $\loc(\hat t) \in g(T_\AA \times \hat T_\AA)$, and commutativity of this diagram, imply that 
$$\sigma(\hat t) \in \Ker \left( \Sha^1(\alg{R}) \To \Sha(\alg{T} \times \alg{\hat T}) \right).$$  
This gives a new characterization of $H^\dag$,
\begin{align*}
H^\dag &= \left\{ t \in T_\AA : \delta j \Kum t \in \Im \left( \hat T_F \xrightarrow{\loc \Kum} P^1(\alg{\hat T}_{[n]}) \right) \right\}, \\
&= \left\{ t \in T_\AA : \delta j \Kum t \in \Im \left( \sigma^{-1}(\Sha^1(\alg{R})) \xrightarrow{\loc \Kum} P^1(\alg{\hat T}_{[n]}) \right) \right\}.
\end{align*}
Consider the subgroup of $H^\dag$:
$$J^\dag =  \left\{ t \in T_\AA : \delta j \Kum t \in \Im \left( \Ker(\sigma) \xrightarrow{\loc \Kum} P^1(\alg{\hat T}_{[n]}) \right) \right\}.$$ 

We find an injective homomorphism and a surjective homomorphism as below.
$$ \frac{H^\dag}{J^\dag} \xhookrightarrow{\delta j \Kum}  \frac{ \Im(\sigma^{-1}(\Sha^1(\alg{R})) \xrightarrow{\loc \Kum}  P^1(\alg{\hat T}_{[n]}))}{\Im(\Ker(\sigma) \xrightarrow{\loc \Kum} P^1(\alg{\hat T}_{[n]}))} \xtwoheadleftarrow{\loc \Kum} \frac{ \sigma^{-1} \Sha^1(\alg{R}) }{\Ker(\sigma)}.$$

We arrive at an estimate of the index,
\begin{equation}
\label{EstimateHJ}
\# (H^\dag / J^\dag) \leq \# \Ker \left( \Sha^1(\alg{R}) \To \Sha^1(\alg{T} \times \alg{\hat T}) \right).
\end{equation}

Finally, consider an element $t \in J^\dag$.  We find that there exists an element $\hat t \in T_F$ such that $\sigma(\hat t) = 0$ and $\delta j \Kum t = \loc \Kum(\hat t)$.  Since $\sigma(\hat t) = 0$, there exist elements $u \in T_F$, $\hat u \in \hat T_F$, such that
$$\hat t = g(u,\hat u) = \delta j(u) \cdot \hat u^{-n}.$$
Therefore
$$\delta j \Kum t = \loc \Kum \left( \delta j(u) \cdot \hat u^{-n} \right).$$
Since $\delta j \Kum = \Kum \delta j$, and $\loc \delta j = \delta j \loc$, we find
$$\delta j \Kum \left( t \cdot u^{-1} \right) = \Kum \loc \hat u^{-n}.$$
(Here we write $u$ instead of $\loc(u)$ to identify an element of $T_F$ with its image in $T_\AA$.)  Hence, for all $s \in T_\AA$,
$$\Comm ( t u^{-1}, s) = \TD \left( \delta j \Kum (t u^{-1} ), \Kum s \right) = \TD \left(  \Kum \loc \hat u, \Kum s \right)^{-n} = 1.$$
Therefore $t u^{-1} \in Z_\AA^\dag$ and $t \in T_F \cdot Z_\AA^\dag$.  We have demonstrated that $J^\dag \subset T_F Z_\AA^\dag$.  

The diagram below summarizes the containments among subgroups of $T_\AA$ and their indices from \eqref{EstimateNE}, \eqref{EstimateEH}, \eqref{EstimateHJ}.
\begin{center}
\begin{tikzpicture}
\node (TZ) at (0,0) {$T_F Z_\AA^\dag$};
\node (N) at (2,2.25) {$N^\dag$};
\node (E) at (2,0.75) {$E^\dag Z_\AA^\dag$};
\node (H) at (2,-0.75) {$H^\dag Z_\AA^\dag$};
\node (J) at (2,-2.25) {$J^\dag Z_\AA^\dag$};
\draw (TZ) -- (N);
\draw (TZ) -- (J);
\draw (N) to node[right] {$\leq \#  \left( \Sha^1(F, \alg{T})_{[n]} / \Im(\Sha^1(\alg{T}_{[n]})) \right)$} (E);
\draw (E) to node[right] {$\leq \#  \left( \Sha^1(F, \alg{\hat T})_{[n]} / \Im(\Sha^1(\alg{\hat T}_{[n]})) \right)$} (H);
\draw (H) to node[right] {$\leq \# \Ker \left( \Sha^1(\alg{R}) \To \Sha^1(\alg{T} \times \alg{\hat T}) \right)$} (J);
\end{tikzpicture}
\end{center}

These estimates demonstrate, by a squeeze, that $\tilde N = Z_{\tilde T_\AA}(T_F)$ contains the abelian subgroup $T_F \tilde Z_\AA$ with finite index (bounded by the product of the three finite indices above).  Therefore (reducing the strength of our estimates a bit), there exists a maximal abelian subgroup $\tilde M \subset \tilde T_\AA$ such that 
$$T_F \tilde Z_\AA \subset \tilde M \subset \tilde N, \text{ and }$$
$$\# \frac{\tilde M}{T_F \tilde Z_\AA} \leq \# \Ker \left( \Sha^1(\alg{R}) \To \Sha^1(\alg{T} \times \alg{\hat T}) \right) \cdot \# \Sha^1(\alg{T} \times \alg{\hat T})_{[n]}.$$
If $\Sha^1(\alg{T})_{[n]} = \Sha^1(\alg{\hat T})_{[n]} = \Sha^1(\alg{R}) = 0$, then $T_F \tilde Z_\AA$ is a maximal abelian subgroup of $\tilde T_\AA$.

In particular, if $\alg{T}$ is a split torus, then $\alg{\hat T}$ also is split, and Hilbert's Theorem 90 implies that $H_{\et}^1(F, \alg{T}) = H_{\et}^1(F, \alg{\hat T}) = 0$.  Thus $\Sha^1(\alg{T}) = \Sha^1(\alg{\hat T}) = 0$.  When $\alg{T}$ and $\alg{\hat T}$ are split, $\alg{R}$ is a split group scheme of multiplicative type.  The component group of $\alg{R}$ is $\alg{\hat \nu}$, a finite $n$-torsion group.  The structure theorem for finitely-generated abelian groups (applied to the character lattice of $\alg{R}$) implies that there exists an isomorphism of group schemes over $F$:
$$\alg{R} \isom \alg{G}_m^r \times \prod_{i=1}^s \alg{\mu}_{n_i},$$
where $n_1, \ldots, n_s$ are divisors of $n$.  Since $F$ contains the $n^{\th}$ roots of unity, the Grunwald-Wang theorem implies that $\Sha^1(\alg{\mu}_{n_i}) = 0$ for all $i$.  Hilbert's Theorem 90 implies that $\Sha^1(\alg{G}_m) \subset H^1(F, \alg{G}_m) = 0$, and so $\Sha^1(\alg{R}) = 0$.
\qed

\begin{corollary}
\label{Mult1SplitTori}
Let $\chi$ be a unitary ACC as before.  The representation $L_\chi^2(T_F \backslash \tilde T_\AA)$ decomposes as a finite isotypic direct sum of unitary irreps:
$$\left[ L_\chi^2(T_F \backslash \tilde T_\AA) \right] = m_\chi \cdot [\pi_\chi],$$
where the multiplicity satisfies:
$$m_\chi \leq \# \Ker \left( \Sha^1(\alg{R}) \To \Sha^1(\alg{T} \times \alg{\hat T}) \right) \cdot \# \Sha^1(\alg{T} \times \alg{\hat T})_{[n]}.$$
In particular, if $\alg{T}$ is a split torus, then $m_\chi = 1$.
\end{corollary}
\proof
This is a direct result of the previous Theorem and Proposition \ref{MultiplicityIndex}.
\qed

\subsection{Global pouches}

Let $\Irr_\epsilon^{\unit}(\tilde T_\AA)$ denote the set of unitary $\epsilon$-genuine automorphic representations of $\tilde T_\AA$; these are the irreducible unitary representations of $\tilde T_\AA$ which are equivalent to a summand of $L_\chi^2(T_F \backslash T_\AA)$ for some $\epsilon$-genuine unitary ACC $\chi$.  Such a representation of $\tilde T_\AA$ is determined by its automorphic central character, and so this gives a bijection
\begin{equation}
\label{PI}
\Irr_\epsilon(\alg{T}) \leftrightarrow \Hom_\epsilon(T_F \cap \tilde Z_\AA \backslash \tilde Z_\AA, U(1)), \quad [\pi_\chi] \mapsto \chi.
\end{equation}

Recall that $\alg{T}^\sharp$ is the torus with character lattice $\sheaf{Y}^\sharp$, fitting into a short exact sequence of groups over $F$, $\alg{\mu} \Into \alg{T}^\sharp \xtwoheadrightarrow{i} \alg{T}$.  Write $\tilde T_\AA^\sharp$ for the extension of $T_\AA^\sharp$ obtained by pulling back $\tilde T_\AA$.  Locally, we have $i(T_v^\sharp) \subset Z_v^\dag$, and so we have a homomorphism
$$i_{\AA / F} \From \frac{T_\AA^\sharp}{T_F^\sharp} \To \frac{Z_\AA^\dag}{Z_\AA^\dag \cap T_F}.$$
\begin{definition}
Two unitary ACCs $\chi_1$ and $\chi_2$ belong to the same \defined{pouch} if their pullbacks to $T_F^\sharp \backslash \tilde T_\AA^\sharp$ coincide.  Two automorphic representations $\pi_1, \pi_2 \in \Irr_\epsilon^{\unit}(\tilde T_\AA)$ belong to the same \defined{pouch} if their central characters belong to the same pouch.
\end{definition}

Observe that if $\chi_1$ and $\chi_2$ belong to the same pouch, then their local components $\chi_{1,v}$ and $\chi_{2,v}$ pull back to the same characters of $\tilde T_v^\sharp$ for all places $v \in \VV$.  Hence if $\pi_1, \pi_2 \in \Irr_\epsilon^{\unit}(\tilde T_\AA)$ belong to the same pouch, then their local factors $\pi_{1,v}$, $\pi_{2,v}$ belong to the same local pouches.

A careful computation may lead to a precise description of global pouches, but this is left for another article.  When $\alg{T}$ is split, we can say something precise.
\begin{proposition}
Suppose that $\alg{T}$ is split.  Then every pouch is a singleton, and pulling back the central character via $i$ parameterizes the unitary automorphic representations of $\tilde T_\AA$:
$$ \Irr_\epsilon^{\unit}(\tilde T_\AA) \leftrightarrow \Hom_\epsilon \left( \frac{\tilde T_\AA^\sharp }{\alg{\mu}(\AA) \cdot T_F^\sharp}, U(1) \right), \quad [\pi_\chi] \mapsto \chi \circ i_{\AA / F}.$$
\end{proposition}
\proof
The short exact sequence of commutative groups over $F$, $\alg{\mu} \Into \alg{T}^\sharp \xtwoheadrightarrow{i} \alg{T}$, yields a commutative diagram with exact rows:
$$\begin{tikzcd}
\alg{\mu}(F) \inarrow{r} \inarrow{d} & T_F^\sharp \arrow{r}{i} \arrow{d}{\loc} & T_F \arrow{r}{\partial} \arrow{d}{\loc} & H_{\et}^1(F, \alg{\mu}) \arrow{d}{\loc} \\
\alg{\mu}(\AA) \inarrow{r} & T_\AA^\sharp \arrow{r}{i} & T_\AA \arrow{r} \arrow{r}{\partial} & P^1(\alg{\mu})
\end{tikzcd}$$ 

If $z \in Z_\AA^\dag$, then $z \in \Im(T_\AA^\sharp \To T_\AA)$ by Theorem \ref{GlobalCenterImage}.  It follows that the homomomorphism $i_{\AA / F} \From T_\AA^\sharp / T_F^\sharp \To Z_\AA^\dag / (Z_\AA^\dag \cap T_F)$ is surjective.  

Next, suppose that $z \in Z_\AA^\dag \cap T_F$.  The commutative diagram demonstrates that $\loc(z) \in \Im(T_\AA^\sharp \To T_\AA)$ and a diagram chase demonstrates that $\partial z \in \Sha^1(\alg{\mu})$.  Since $\alg{T}$ is split, $\alg{\mu}$ is a split $n$-torsion group, and we have $\Sha^1(\alg{\mu}) = 0$ by the Grunwald-Wang theorem.  Thus $\partial z = 0$ and so $z \in \Im(T_F^\sharp \To T_F)$.  We have demonstrated that
$$i(T_F^\sharp) = Z_\AA^\dag \cap T_F.$$

Finally, suppose that $ t \in T_\AA^\sharp$ and $i(t) \in Z_\AA^\dag \cap T_F$.  Let $t' \in T_F^\sharp$ be an element for which $i(t') = i(t)$, and so $i(t / t') = 1 \in Z_\AA^\dag$.  It follows that
$$t / t' \in \alg{\mu}(\AA).$$
We end up with a short exact sequence
$$\frac{ \alg{\mu}(\AA)}{\alg{\mu}(F)} \Into \frac{ T_\AA^\sharp}{T_F^\sharp} \xtwoheadrightarrow{i_{\AA / F}} \frac{ Z_\AA^\dag}{Z_\AA^\dag \cap T_F}.$$
The proposition follows by Pontrjagin duality and the bijection \eqref{PI}.
\qed

For split tori over local fields, we found an embedding in Theorem \ref{ParamLocalTori},
$$\Irr_\epsilon(\tilde T_v) \hookrightarrow \Irr_\epsilon(\tilde T_v^\sharp).$$
Now for split tori over global fields, we find an embedding,
$$\Irr_\epsilon^{\unit}(\tilde T_\AA) \hookrightarrow \Irr_\epsilon^{\unit}(\tilde T_\AA^\sharp).$$
Both are obtained by pulling back via $i \From \alg{T}^\sharp \To \alg{T}$, and are thus compatible with factorization of automorphic representations.  For nonsplit tori, we obtain not an embedding, but a finite-to-one map whose nonempty fibres are local or global pouches.  We hope to examine the fine structure of these pouches in a later paper, as we believe they should behave analogously to (or perhaps as a refinement of) L-packets for covers of tori.

\bibliography{CovLang2014.bib}

\end{document}